\documentclass[12pt,leqno]{article}
\usepackage{graphicx}
\usepackage{hyperref}
\usepackage{amssymb,amsmath,amsthm}
\usepackage{epstopdf}
\usepackage{mathrsfs}
\usepackage{xypic}
\usepackage{comment}
\usepackage{epsfig}
\usepackage{color}
\usepackage{tabls}

\usepackage{float}

\theoremstyle{plain}
\newtheorem{theorem}{Theorem}[section]
\newtheorem{conjecture}[theorem]{Conjecture}

\newtheorem{lemma}[theorem]{Lemma}
\newtheorem{proposition}[theorem]{Proposition}

\theoremstyle{definition}
\newtheorem{definition}[theorem]{Definition}

\usepackage{cancel}

\def\FF{\mathbb{F}}
\def\NN{\mathbb{N}}
\def\RR{\mathbb{R}}
\def\QQ{\mathbb{Q}}
\def\ZZ{\mathbb{Z}}

\def\L{\mathcal{L}}
\def\R{\mathcal{R}}

\def\Tsf{\mathsf{T}}
\def\Jsf{\mathsf{J}}
\def\maj{\mathsf{maj}}
\def\DR{\mathsf{DR}}

\def\Box{\operatorname{Box}}
\def\Cat{\operatorname{Cat}}
\def\hull{\operatorname{hull}}
\def\area{\operatorname{area}}
\def\bounce{\operatorname{bounce}}
\def\dinv{\operatorname{dinv}}
\def\sweep{\operatorname{sweep}}

\def\Dyck{\operatorname{Dyck}}

\newcommand{\sqbinom}[2]{\genfrac{[}{]}{0pt}{}{#1}{#2}}

\def\a{\mathbf{a}}

\def\x{\mathbf{x}}
\def\y{\mathbf{y}}
\def\z{\mathbf{z}}
\def\0{\mathbf{0}}

\def\tbinom#1#2{\bigl(\hskip-1.5pt\begin{smallmatrix}#1\\#2\end{smallmatrix}\hskip-1.5pt\bigr)}
\def\tsqbinom#1#2{\bigl[\hskip-1.5pt\begin{smallmatrix}#1\\#2\end{smallmatrix}\hskip-1.5pt\bigr]}
\def\pmod#1{\ (\textup{mod }#1)}
\def\tmod#1{\ \textup{mod }#1}
\def\downstrut{\noindent\lower5pt\vbox{}}
\def\upstrut{\noindent \vbox to 12pt{}}

\title{Lattice Points and Rational $q$-Catalan Numbers}
\author{Drew Armstrong}
\date{September 2025}

\allowdisplaybreaks

\begin{document}

\maketitle

\begin{abstract}
For each pair of coprime integers $a$ and $b$ we have the {\em rational $q$-Catalan number}  $\Cat(a,b)_q=\bigl[\hskip-1.5pt \begin{smallmatrix}{a+b}\\{a}\end{smallmatrix}\hskip-1pt\bigr]_q/[a+b]_q$. It is known that this is a polynomial in $q$ with non-negative integer coefficients, but the nature of these coefficients is still mysterious. Our current understanding is based on the {\em rational shuffle conjecture} of Bergeron, Garsia, Leven and Xin from 2014, which was proved by Mellit in 2016, building on earlier work with Carlsson. This theorem realizes $\Cat(a,b)_q$ as the generating function for the statistic ``$\area-\dinv+(a-1)(b-1)/2$'' defined on rational Dyck paths. However, this statistic is difficult to work with and leaves some phenomena unexplained. For example, it does not prove the conjecture that the difference $\Cat(a,c)_q-\Cat(a,b)_q$ has non-negative coefficients whenever $\mathrm{gcd}(a,b)=\mathrm{gcd}(a,c)=1$ and $b<c$. The current paper proposes to look at lattice points instead of Dyck paths. Our idea is to fix $a$ and express everything in terms of the weight lattice $\L$ and root lattice $\R$ of type $A_{a-1}$. Based on ideas of Paul Johnson, we conjecture the existence of certain ``Johnson statistics'' $\Jsf:\R\to\ZZ$ and we prove this conjecture for $a\leq 20$. We show that these statistics satisfy many remarkable properties including a $q$-analogue of Brion's theorem for simplices.
\end{abstract}

{
\setlength{\parskip}{0pt}
\tableofcontents
}

\bigskip

{\em Note: This new version of the paper is substantially revised. It includes new results, new conjectures, new figures and a new numbering system. The most important changes are: (1) The number of axioms for Johnson statistics is reduced from three to two, since we found that one axiom was implied by the other two. (2) The main conjecture is verified for many more cases. (3) We prove an analogue of Brion's theorem for Johnson statistics. (4) We generalize the monotonicity conjecture to the case of non-coprime parameters.}

\section{Introduction} 
\begin{figure}[h]
\begin{center}
\includegraphics{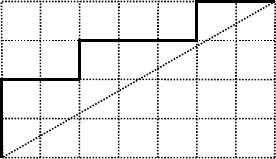}
\end{center}
\caption{A $(4,7)$-Dyck path.}
\label{fig:dyck_path}
\end{figure}

Let $a,b\geq 1$ be integers. An {\em
$(a,b)$-Dyck path} is a lattice path from the point $(0,0)$ to the
point $(b,a)$ that stays weakly above the diagonal.\footnote{We
use $(b,a)$ instead of $(a,b)$ because we want the pictures to be
horizontal when $b>a$.} To be precise, this is a sequence of integer
points $(x_i,y_i)$ for $0\leq i\leq a+b$ such that $(x_0,y_0)=(0,0)$,
$(x_{a+b},y_{a+b})=(b,a)$, $(x_{i+1}-x_i,y_{i+1}-y_i)\in\{
(1,0),(0,1)\}$ and $y_i/x_i\geq b/a$ for all $i$. Figure \ref{fig:dyck_path} shows an example of a $(4,7)$-Dyck path.

Let $\Dyck_{a,b}$ be the set of $(a,b)$-Dyck paths. The case of
coprime $a$ and $b$ is special. It has been known since Bizley
\cite{bizley} that
\begin{equation*}
\#\Dyck_{a,b} = \frac{1}{a+b} \binom{a+b}{a} \quad \text{when
$\gcd(a,b)=1$.}
\end{equation*}
The proof interprets each lattice path in the rectangle as a word
of length $a+b$ with $a$ copies of the letter $u$ (up) and $b$
copies of the letter $r$ (right). For example, the Dyck path in Figure \ref{fig:dyck_path}
corresponds to the word $uurrurrrurr$. The number of such words is
clearly $\tbinom{a+b}{a}$. Bizley considered the action of rotation
on the set of words. Since $\gcd(a,b)=1$, each orbit has size exactly
$a+b$, and Bizley observed that each orbit contains a unique Dyck path.

Bizley also gave a formula for $\#\Dyck_{a,b}$ when $\gcd(a,b)\neq
1$, which is more complicated. However, in the simple case $a=b$ we
have a bijection 
$\Dyck_{a,a}\to \Dyck_{a,a+1}$ defined by appending
a right step $(1,0)$ to the end of the path. Hence
\begin{equation*}
\#\Dyck_{a,a} = \#\Dyck_{a,a+1} = \frac{1}{a+a+1}\binom{a+a+1}{a}
= \frac{1}{2a+1}\binom{2a+1}{a},
\end{equation*}
Note that $\Dyck_{a,a}$ is the set of classical Dyck paths and $\binom{2a+1}{a}/(2a+1)$ is the classical Catalan number. Based on this example we define the notation
\begin{equation*}
\Cat(a,b):= \frac{1}{a+b} \binom{a+b}{a} \quad\text{when
$\gcd(a,b)=1$},
\end{equation*}
and we call this a {\em rational Catalan number}.\footnote{Of
course the number $\Cat(a,b)$ is an integer. It is the slope
$b/a$ that is rational.} Note that we have $\Cat(a,b)=\Cat(b,a)$
because switching up and right steps gives a bijection
$\Dyck_{a,b}\leftrightarrow\Dyck_{b,a}$. However, there is another common formula that obscures the symmetry:
\begin{equation*}
\Cat(a,b)= \frac{1}{a} \binom{a-1+b}{a-1}.
\end{equation*}
This paper is concerned with the geometric meaning behind this expression,
which was originally studied by Haiman~\cite{haiman} in the greater generality of Weyl groups (see
Section~\ref{sec:weyl}). For the purpose of this introduction we give a
quick summary. Let~$\Delta$ be the simplex in $\RR^{a-1}$ with $a$ vertices
$(0,\ldots,0)$,\, $(1,0,\ldots,0)$, \ldots, $(0,\ldots,0,1)$. For
any integer $b\geq 0$ we have the following formula for the number
of integer points in the dilated simplex $b\Delta$:
\begin{equation*}
\#(\ZZ^{a-1}\cap b\Delta) = \binom{a-1+b}{a-1}.
\end{equation*}
We also consider a special sublattice of $\ZZ^{a-1}$ called the {\em
root lattice}:
\begin{equation*}
\R = \{(x_1,\ldots,x_{a-1})\in \ZZ^{a-1}: x_1+2x_2+\cdots+(a-1)x_{a-1}
\equiv 0  \mod a\}.
\end{equation*}
If $\gcd(a,b)=1$ then there is a natural action of $\ZZ/a\ZZ$ on the
set of points $\ZZ^{a-1}\cap b\Delta$, each of whose orbits has size
$a$. Furthermore, each orbit contains a unique point of $\R$; hence
\begin{equation*}
\#(\R\cap b\Delta) =\frac{1}{a} \binom{a-1+b}{a-1}.
\end{equation*}
See Section~\ref{sec:rational_cat} for the details. The $q$-binomial coefficient, as usual, is defined as
\begin{equation*}
\sqbinom{a-1+b}{a-1}_q := \frac{[a-1+b]_q!}{[a-1]_q![b]_q!},
\end{equation*}
where $[n]_q=(1-q^n)/(1-q)$, $[n]_q!=\prod_{i=1}^n [i]_q$ and $[0]_q!=1$. This polynomial can be  interpreted combinatorially as the generating function for
the ``area'' statistic on lattice paths in an $(a-1)\times b$
rectangle. Equivalently, we can interpret $\tsqbinom{a-1+b}{a-1}_q$
as the generating function for a certain statistic on the lattice
$\ZZ^{a-1}$, which we call the {\em tilted height}:
\begin{equation*}
\Tsf(x_1,\ldots,x_{a-1}) = x_1+2x_2+\cdots +(a-1)x_{a-1}.
\end{equation*}
To be precise, in Section~\ref{sec:qbinomial} we observe for all integers $b\geq 0$ that
\begin{equation*}
\sum_{\x\in (\ZZ^{a-1}\cap b\Delta)} q^{\Tsf(\x)} =
\sqbinom{a-1+b}{a-1}_q.
\end{equation*}
When $\gcd(a,b)=1$ it is known that $\tsqbinom{a-1+b}{a-1}_q=[a]_q f(q)$
for some polynomial $f(q)\in\NN[q]$, which we call a {\em rational
$q$-Catalan number}:
\begin{equation*}
\Cat(a,b)_q = \frac{1}{[a]_q} \sqbinom{a-1+b}{a-1}_q \in\NN[q]\quad\text{when
$\gcd(a,b)=1$.}
\end{equation*}
This is not easy to prove, and we do not have a good understanding of the coefficients. Our current best understanding is based on a mixture of two statistics on rational Dyck paths, called area and dinv. The following formula is just one consequence of the {\em Rational Shuffle Theorem} of Mellit \cite{mellit} (see Section~\ref{sec:rational_qcat} for details):
\begin{equation*}
\Cat(a,b)_q = \sum_{P\in \Dyck_{a,b}} q^{\text{area}(P)-\text{dinv}(P)+(a-1)(b-1)/2}.
\end{equation*}
Unfortunately, we find the statistic ``area-dinv'' awkward to work with. For example, it does not explain the following phenomenon, which we state as a conjecture since we have not seen it written down. In Section \ref{sec:greedy} we will prove this conjecture for the infinite set of triples $a,b,c$ with $a\leq 20$. See Conjecture \ref{conj:monotone2} for a generalization to the case $\gcd(a,b)=\gcd(a,c)> 1$.

\begin{conjecture}\label{conj:monotone}
For any $a,b,c\geq 1$ with $\gcd(a,b)=\gcd(a,c)=1$ and $b<c$ we have
\begin{equation*}
\Cat(a,c)_q - \Cat(a,b)_q \in\NN[q].
\end{equation*}
\end{conjecture}
For example,
\begin{align*}
\Cat(5,3)_q - \Cat(5,2)_q &= (q^8 + q^6 + q^5 + q^4 + q^3 + q^2 + 1)-(q^4 + q^2 + 1) \\
&= q^8 + q^6 + q^5 + q^3 \in\NN[q].
\end{align*}

Conjecture \ref{conj:monotone} is the weakest conjecture in this paper. The stronger conjectures are based on the following framework, which we hope will lead to a simpler understanding of rational $q$-Catalan numbers. We define the {\em tilted partial order} on the lattice  $\ZZ^{a-1}$ by setting
\begin{equation*}
(x_1,\ldots,x_{a-1}) \leq (y_1,\ldots,y_{a-1}) \quad\Longleftrightarrow\quad x_i+\cdots+x_{a-1}\leq y_i+\cdots+y_{a-1} \text{ for all $1\leq i\leq a-1$}.
\end{equation*}
This poset is graded with rank function $\Tsf$. Let $S=\{\x_0<\cdots<\x_{a-1}\}$ be a chain in this poset that is saturated in the sense that $T(\x_i)=T(\x_0)+i$ for all $0\leq i\leq a-1$. If $k:= T(\x_0)$ then we note that
\begin{equation*}
\sum_{i=0}^{a-1} q^{T(\x_i)} = q^k +q^{k+1} + \cdots + q^{k+a-1} = q^k [a]_q.
\end{equation*}
Any saturated chain in $\ZZ^{a-1}$ of length $a$ is called a {\em ribbon}.\footnote{We borrow this notation from the theory of ribbon tableaux  \cite{llt}. See Figure  \ref{fig:box_22} for the connection.} If $\gcd(a,b)=1$ then we know that the number of points $\#(\ZZ^{a-1}\cap b\Delta)=\binom{a-1+b}{a-1}$ is divisible by $a$. Suppose we have a partition of the subposet $\ZZ^{a-1}\cap b\Delta$ into
$\Cat(a,b)=\tbinom{a-1+b}{a-1}/a$ ribbons. If $\{S_i\}$ is any such partition and
$q^{k_i}[a]_q$ is the sum of $q^{\Tsf(\x)}$ over the points $\x\in S_i$ then we obtain
\begin{align*}
\Cat(a,b)_q &= \frac{1}{[a]_q} \sqbinom{a-1+b}{a-1}_q
= \frac{1}{[a]_q} \sum_{\x\in (\ZZ^{a-1}\cap b\Delta)} q^{\Tsf(\x)}\\
&= \frac{1}{[a]_q} \sum_{i=1}^{\Cat(a,b)}\sum_{\x\in S_i}
q^{\Tsf(\x)}
= \frac{1}{[a]_q} \sum_{i=1}^{\Cat(a,b)} q^{k_i}[a]_q
= \sum_{i=1}^{\Cat(a,b)} q^{k_i}.
\end{align*}
Thus each ribbon partition of the poset $\ZZ^{a-1}\cap
b\Delta$ gives a combinatorial interpretation of $\Cat(a,b)_q$. We conjecture that at least one ribbon partition always exists. More specifically, Conjecture \ref{conj:ribbon} implies that there exists at least one ribbon partition $\mathcal{S}=\{S\}$ of the infinite lattice
$\ZZ^{a-1}$ with the following properties:
\begin{itemize}
\item (Periodic) If $S\in\mathcal{S}$ then for any
translation $\y\in\ZZ^{a-1}$ we have $(S+a\y) \in\mathcal{S}$.
\item (Number Theoretic) If $\gcd(a,b)=1$ then the set $\ZZ^{a-1}\cap
b\Delta$ is a union of ribbons.
\end{itemize}
Together these properties imply that no ribbon of the partition can cross one of the hyperplanes $x_i=ak+\varepsilon$ with $k\in\ZZ$ or $x_1+\cdots+x_{a-1}=b+\varepsilon$, where $\gcd(a,b)=1$ and $0<\varepsilon<1$. The existence of such a ribbon partition naturally implies Conjecture~\ref{conj:monotone} since it implies that the set $\ZZ^{a-1}\cap (c\Delta-b\Delta)$ is a union of ribbons.

Our concept of ribbon partitions was inspired by work of Paul Johnson \cite{johnson}. Recall that the root lattice $\R$ is the set of points $\x\in\ZZ^{a-1}$ with $\Tsf(\x)\equiv 0$ mod $a$. Each ribbon necessarily contains a unique point in $\R$; hence any ribbon partition of $\ZZ^{a-1}$ induces a statistic $\Jsf:\R\to\ZZ$ defined by setting $\Jsf(\x):=\Tsf(\y)$, where $\y$ is the ``lowest''
point in the ribbon containing $\x$. If our ribbon partition satisfies
the above two properties then we will see that the statistic $\Jsf$ satisfies the
corresponding properties:
\begin{itemize}
\item (Periodic) For any $\x\in \R$ and $\y\in\ZZ^{a-1}$ we have
$\Jsf(\x+a\y)=\Jsf(\x)+a\Tsf(\y)$.
\item (Catalan) For any integer $b\geq 0$ coprime to $a$ we have
\begin{equation*}
\Cat(a,b)_q = \sum_{\x \in (\R\cap b\Delta)} q^{\Jsf(\x)}.
\end{equation*}
\end{itemize}
The existence of such a function $\Jsf:\R\to\ZZ$ is closely related
to a conjecture of Johnson (see Section~\ref{sec:johnson}), so we call
it a {\em Johnson statistic}. In Section~\ref{sec:ribbons} we will see that there exists a unique
Johnson statistic when $a=3$ and there is a geometrically nice
Johnson statistic when $a=4$. We prove by computation that ribbon partitions exist for $a=5$, but their structure is not so clear.

In Section \ref{sec:greedy} we discuss a weaker conjecture in which we replace ribbons by {\em standard sets}, i.e., sets of $a$ elements with consecutive tilted heights that need not be chains in the tilted partial order. By using a greedy algorithm we can prove that such a ``standard partition'' exists for all $a\leq 20$. This weaker result is still enough to prove the existence of Johnson statistics, and hence Conjecture \ref{conj:monotone}, for all $a\leq 20$.

These concepts can also be framed in terms of Young diagrams. In Section~\ref{sec:qbinomial} we observe that the tilted partial order on the set
$\ZZ^{a-1}\cap b\Delta$ is isomorphic to the subposet of Young's
lattice that Stanley called $L(a-1,b)$ \cite{stanley_paper}, whose rank numbers are the
coefficients of $\tsqbinom{a-1+b}{a-1}_q$. These posets are known
to be extremely subtle. Take, for example, the fact that the rank
numbers are unimodal. Sylvester~\cite{sylvester} gave the first proof
of this, calling it ``a theorem that has been awaiting proof for the
last quarter of a century and upwards''. The most elementary proof is
due to Proctor~\cite{proctor}, which is a simplification of Stanley's
proof from \cite{stanley_paper}. The ideal combinatorial proof would
be to construct a symmetric chain decomposition (SCD) of the poset
$L(a-1,b)$, i.e., a partition of the poset into saturated chains,
each of which is centered at the middle rank. Despite much effort
since the 1980s this problem is still open.

We view the monotonicity property of rational $q$-Catalan numbers (Conjecture \ref{conj:monotone}) as roughly analogous to the unimodality property of $q$-binomial coefficients, and the existence of ribbon partitions (Conjecture \ref{conj:ribbon}) as roughly analogous to the existence of symmetric chain decompositions. Due to this superficial similarly, it could be that these conjectures are difficult.

\section{Root and weight lattice notation}
\label{sec:lattices}

In this section we establish notation for the root system of type $A$. See \cite{humphreys} for the general theory of root systems. Throughout the paper we fix an integer $a\geq 2$ and we
work within the {\em weight lattice of type $A_{a-1}$}. We denote
this lattice by $\L=\ZZ^{a-1}$, and we express everything in terms of the basis of  {\em fundamental weights}
 \begin{equation*}
 \omega_1=(1,0,\ldots,0),\quad \ldots,\quad \omega_{a-1}=(0,\ldots,0,1),
 \end{equation*}
 which we view as column vectors. We equip $\L$ with
 the inner product $\langle \x,\y\rangle = \x^TC^{-1}\y$, where $C$
 is the $(a-1)\times (a-1)$ {\em Cartan matrix}
 {\small
  \begin{equation*}
  \def\pmin{\phantom{-}}
 C := \begin{pmatrix} \pmin 2 & -1 & \\ -1 & \pmin 2 & -1 & & \\ & -1 & \pmin 2 & -1 &
 & \\ & &  \ddots &  \ddots &  \ddots  \\  &&& -1 & \pmin 2 & -1\\ &&&& -1 &
 \pmin 2\end{pmatrix}.
  \end{equation*}}%
 This inner product is typically called the {\em Killing
 form}. Let $\alpha_i$ denote the $i$-th column vector of $C$,
 so that $\langle \omega_i,\alpha_j\rangle =\delta_{i,j}$. The vectors
 $\alpha_1,\ldots,\alpha_{a-1}\in \L$ are called the {\em simple roots},
 and they define a sublattice of $\L$ called the {\em root lattice}:
 \begin{equation*}
 \R:=\ZZ\alpha_1+\cdots+\ZZ\alpha_{a-1} \subseteq \L.
 \end{equation*}
The root lattice has a convenient interpretation in terms of weight
coordinates:
 \begin{equation*}
 \R=\{(x_1,\ldots,x_{a-1})\in\ZZ^{a-1}:
 x_1+2x_2+\cdots+(a-1)x_{a-1}\equiv 0 \mod a\}.
 \end{equation*}
More generally, for each integer $0\leq k\leq a-1$ we define the set
$\R_k\subseteq \L$ by
  \begin{equation*}
 \R_k=\{(x_1,\ldots,x_{a-1})\in\ZZ^{a-1}:
 x_1+2x_2+\cdots+(a-1)x_{a-1}\equiv k \mod a\},
 \end{equation*}
 so that $\R_0=\R$. Note that $\L$ is partitioned by the sets
 $\R_0,\ldots,\R_{a-1}$. In fact, these are the cosets of the subgroup
 $\R\subseteq \L$, and the map $k\mapsto \R_k$ is a well-defined
 group isomorphism $\ZZ/a\ZZ\cong \L/\R$. It follows from this that $\det(C)=a$.
 
We will also need the following polytopes, called the {\em positive cone} $K$, the {\em fundamental alcove} $\Delta$ and the {\em fundamental parallellotope} $\Pi$. These are defined as polytopes living in the real vector space $\L\otimes\RR = \RR^{a-1}$ spanned by the weight lattice:
\begin{align*}
K &:= \{(x_1,\ldots,x_{a-1}): 0\leq x_i \text{ for all $i$\,}\},\\
\Delta &:= \{(x_1,\ldots,x_{a-1}): 0\leq x_i \text{ for all $i$ and } \sum x_i\leq 1\},\\
\Pi &:= \{(x_1,\ldots,x_{a-1}): 0\leq x_i< 1 \text{ for all $i$\,}\}.
 \end{align*}
Note that $\Delta$ is the convex hull of the basis vectors $\omega_1,\ldots,\omega_{a-1}$ together with the origin. The set $\Pi$ is a fundamental domain for the translation action of $\L$ on $\L\otimes\RR$. We also define the {\em fundamental box}
\begin{equation*}
\Box:=a\Pi=\{(x_1,\ldots,x_{a-1}): 0\leq x_i< a \text{ for all $i$\,}\},
\end{equation*}
which is a fundamental domain for the translation action of the dilated lattice $a\L$ on $\L\otimes\RR$. In Section \ref{sec:johnson} we will study the sets of lattice points $\L\cap \Box$ and $\R\cap \Box$.
 
\section{The tilted partial order}
\label{sec:tilted}
  
  Based on the description of the
 root lattice in the previous section, we define the {\em tilted
 height function} $\Tsf:\L\to\ZZ$ as 
\begin{equation*}
\Tsf(x_1,\ldots,x_{a-1}):=x_1+2x_2+\cdots+(a-1)x_{a-1}.
\end{equation*}
Note that the coset $\R_k$ of the root lattice consists of all
points with tilted height $\equiv k \tmod a$. Note also that the zero
vector $0\in\L$ satisfies $\Tsf(0)=0$ and the fundamental weight basis
satisfies $\Tsf(\omega_i)=i$ for all $1\leq i\leq a-1$. Alternatively,
we can define the tilted height via the inner product with the last
fundamental weight:
\begin{equation*}
\Tsf(\x) = a\langle \x,\omega_{a-1}\rangle = a \x^T C^{-1}
\omega_{a-1}.
\end{equation*}
Thus $\Tsf(\x)$ measures the height of a point in the direction of $\omega_{a-1}=(0,\ldots,0,1)$ with respect to the Killing form. Next we define a partial order on the set $\L$. Consider the
{\em tilted basis}
\begin{equation*}
\{\omega_1, \omega_2-\omega_1,\ldots, \omega_{a-1}-\omega_{a-2}\},
\end{equation*}
which is unimodular with respect to the standard basis. Given a point $\x=\sum_i x_i\omega_i \in\L$ in standard
coordinates we note that
\begin{equation*}
\x= (x_1+\cdots+x_{a-1})\omega_1 +
(x_2+\cdots+x_{a-1})(\omega_2-\omega_1) + \cdots +
(x_{a-1})(\omega_{a-1}-\omega_{a-2}).
\end{equation*}
Then for all points $\x,\y\in\L$ we define the {\em tilted partial order} as the
componentwise order on coefficients with respect to the tilted basis:
\begin{equation*}
\x\leq \y \qquad\iff \qquad  x_i+\cdots +x_{a-1} \leq
y_i+\cdots+y_{a-1} \quad\text{ for all $1\leq i\leq a-1$.}
\end{equation*}
Observe that $\x\prec \y$ is a cover relation if and only if $\y-\x$
is an element of the tilted basis. Furthermore, observe that each
element of the tilted basis has tilted height $1$. Combining this
fact with the $\ZZ$-linearity of the tilted height function shows that
\begin{equation*}
\x\prec\y \quad\implies\quad \Tsf(\y)=\Tsf(\x)+1.
\end{equation*}
In other words, $\Tsf$ is a rank function for $(\L,\leq)$. Finally,
since the poset $(\L,\leq)$ is isomorphic to the componentwise partial
order on $\ZZ^{a-1}$ we conclude that this partial order has meets
and joins. Thus we can regard $\L$ as a ``lattice'' in the double sense
of being both a discrete additive subgroup of a real vector space,
and also as a poset with meets and joins.

\begin{figure}[h]
\begin{center}
\includegraphics[scale=0.9]{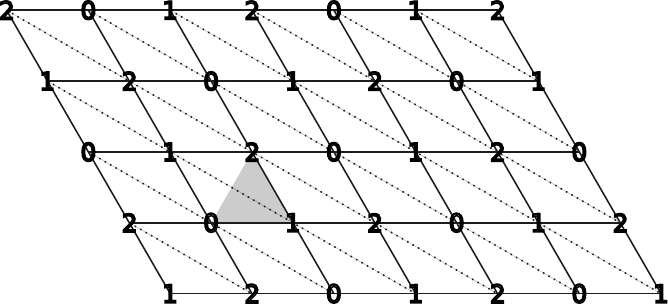}
\end{center}
\caption{The tilted partial order on the weight lattice for $a=3$.}
\label{fig:A2_lattice}
\end{figure}

Figure~\ref{fig:A2_lattice} shows the tilted partial order on $\L$
in the case $a=3$, with its points labeled by their tilted heights
modulo $3$. The vertices labeled by $0$ are the points of the root
lattice. The solid lines are cover relations in the tilted order. The
dotted lines connect points with the same tilted height. The
shaded triangle is the fundamental alcove, which has vertices
$0$, $\omega_1$, $\omega_2$. The angles in the picture are determined by the Killing form. So, for example, the angle $\theta$ between the basis vectors $\omega_1,\omega_2$ satisfies
\begin{equation*}
\cos\theta = \frac{\langle \omega_1,\omega_2\rangle}{\sqrt{\langle \omega_1,\omega_1\rangle }\,\,\sqrt{\langle \omega_2,\omega_2\rangle }} = \frac{1/3}{\sqrt{2/3}\,\, \sqrt{2/3}} = \frac{1}{2}.
\end{equation*}

\section{Rational Catalan numbers}
\label{sec:rational_cat}

Recall that the fundamental alcove is the convex hull of the
origin and the standard basis vectors within the real vector space
spanned by $\L$:
\begin{equation*}
\Delta = \hull\{0,\omega_1,\ldots,\omega_{a-1}\} \subseteq \RR\otimes\L.
\end{equation*}
We are interested in the dilations $b\Delta$ for various integers
$b\geq 0$. The points of $\L$ inside $b\Delta$ have the explicit description
\begin{equation*}
\L\cap b\Delta = \{(x_1,\ldots,x_{a-1})\in\ZZ^{a-1}: 0\leq x_i
\text{ for all $i$ and } \sum_{i=1}^{a-1} x_i \leq b\}.
\end{equation*}
The number of such points is
\begin{equation*}
\#(\L\cap b\Delta) = \binom{a-1+b}{a-1}.
\end{equation*}
We also consider a certain action of $\ZZ/a\ZZ$ on the set $\L\cap
b\Delta$. First we introduce a slack variable $x_0$ and consider the set
\begin{equation*}
X_b=\{(x_0,x_1,\ldots,x_{a-1})\in\ZZ^{a}: 0\leq x_i \text{
for all $i$ and } \sum_{i=0}^{a-1} x_i = b\}.
\end{equation*}
Note that the projection $\pi:X_b\to \L\cap b\Delta$ defined by
$\pi(x_0,x_1,\ldots,x_{a-1}):=(x_1,\ldots,x_{a-1})$ is a bijection
with inverse  $\pi^{-1}(x_1,\ldots,x_{a-1})=(b-\sum_{i=1}^{a-1}
x_i,x_1,\ldots,x_{a-1})$. We define an action of $\ZZ/a\ZZ$ on
the set $\L\cap b\Delta$ by pulling back the obvious action
of $\ZZ/a\ZZ$ on the set $X_b$:
\begin{equation*}
\rho(x_0,x_1,\ldots,x_{a-1}) = (x_{a-1},x_0,\ldots,x_{a-2}).
\end{equation*}
That is, we define the function $\phi_b:\L\cap b\Delta \to \L\cap b\Delta$ by $\phi_b:=\pi \circ \rho \circ \pi^{-1}$, so that\footnote{One can check that the extension $\phi_b:\RR^{a-1}\to \RR^{a-1}$ is an isometry with respect to the Killing form. Thomas and Williams~\cite{thomas} studied the
action of this map on the centroids of the $b^{a-1}$ alcoves contained
in the dilation $b\Delta$.}
\begin{align*}
\phi_b(x_1,\ldots,x_{a-1}) &= (\pi \circ \rho \circ \pi^{-1})(x_1,\ldots,x_{a-1})\\
&= (\pi \circ\rho)(b-\textstyle{\sum_{i=1}^{a-1}
x_i},x_1,\ldots,x_{a-1})\\
&=  \pi (x_{a-1},b-\textstyle{\sum_{i=1}^{a-1} x_i},x_1,\ldots,x_{a-2})\\
&= (b-\textstyle{\sum_{i=1}^{a-1} x_i},x_1,\ldots,x_{a-2}).
\end{align*}
We observe that the map $\phi_b$ shifts the tilted height of a
given point by $b \tmod a$. Indeed, for any point $\x\in\L\cap b\Delta$
we have
\begin{align*}
\Tsf(\phi_b(\x)) &=\Tsf(b-\textstyle{\sum_{k=1}^{a-1} x_k},x_1,\ldots,x_{a-2}) \\
&= (b-\textstyle{\sum_{i=1}^{a-1}x_i})+2x_1+3x_2+\cdots+(a-1)x_{a-2} \\
&= b-x_{a-1} + x_1+2x_2+\cdots+(a-2)x_{a-2} \\
&= b-x_{a-1} + \Tsf(\x) - (a-1)x_{a-1}\\
&= \Tsf(\x)+b-ax_{a-1} \\
&\equiv \Tsf(\x)+b \pmod a.
\end{align*}
If $b$ is coprime to $a$ then this implies that the set $\L\cap
b\Delta$ contains an equal number of points from each coset of the
root lattice:\footnote{See Section \ref{sec:non_coprime} for a discussion of the case $\gcd(a,b)\neq 1$.}
\begin{equation*}
\#(\R_0\cap b\Delta)=\#(\R_1\cap b\Delta)=\cdots = \#(\R_{a-1}\cap
b\Delta).
\end{equation*}
In particular, exactly $1/a$ of these points come from the root lattice itself. In summary, for any coprime integers $\gcd(a,b)=1$ the number of points of $\R$ in the dilated simplex $b\Delta$ is the {\em rational Catalan number}:
\begin{equation*}
\#(\R\cap b\Delta)=\Cat(a,b)=\frac{1}{a}\binom{a-1+b}{a-1}.
\end{equation*}

Rational Catalan numbers have a long history. As mentioned in the
introduction, the earliest known appearance is Bizley~\cite{bizley},
where they appear as the number of $(a,b)$-Dyck paths. They were
later rediscovered by Anderson~\cite{anderson} in the context of
representation theory of the symmetric group, where they count the number
of simultaneous {\em $(a,b)$-core partitions}. The interpretation of
$\Cat(a,b)$ in terms of the root lattice was discovered by Haiman
\cite{haiman}, who placed these numbers in the context of the ring
of diagonal coinvariants of a Weyl group. See the next section for more details.

\section{Rational $q$-Catalan and $q,t$-Catalan numbers}
\label{sec:rational_qcat}

\begin{figure}[h]
\begin{center}
	\includegraphics[scale=0.9]{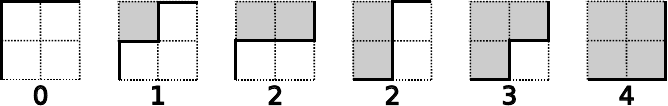}
\end{center}
\caption{Illustration of $\tsqbinom{4}{2}_q=1+q+2q^2+q^3+q^4$.}
\label{fig:gauss_binomial}
\end{figure}

There are several $q$-analogues of Catalan numbers. The
one we have in mind was introduced by F\"urlinger and Hofbauer
\cite{furlinger}. Recall that the $q$-binomial coefficients satisfy the recurrence
\begin{equation*}
\sqbinom{n}{k}_q=q^k\sqbinom{n-1}{k}_q+\sqbinom{n-1}{k-1}_q,
\end{equation*}
from which it follows that $\tsqbinom{n}{k}_q$ is a polynomial in $q$ with non-negative integer coefficients. The nicest combinatorial interpretation of the coefficients is via
lattice paths: the coefficient of $q^\ell$ in $\tsqbinom{n}{k}_q$
is the number of lattice paths in a $k\times (n-k)$ rectangle
having ``area $\ell$\,'', i.e., having $\ell$ unit squares above the
path. Figure~\ref{fig:gauss_binomial} illustrates the expansion
$\tsqbinom{4}{2}_q=1+q+2q^2+q^3+q^4$.

\begin{figure}[h]
\begin{center}
	\includegraphics[scale=0.8]{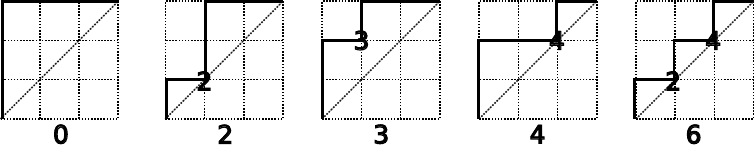}
\end{center}
\caption{Illustration of $\Cat(3,4)_q=1+q^2+q^3+q^4+q^6$.}
\label{fig:major_index}
\end{figure}

The F\"urlinger--Hofbauer $q$-Catalan numbers are defined as 
\begin{equation*}
\Cat(n)_q= \frac{1}{[n]_q}\sqbinom{2n}{n-1}_q.
\end{equation*}
It turns out that this is also a polynomial with non-negative integer
coefficients. The most elementary
proof uses the ``major index of Dyck paths'', which according to
F\"urlinger and Hofbauer \cite{furlinger} goes back to MacMahon and
Aissen. Given an $n\times n$ Dyck path $P\in\Dyck_{n,n}$, the major
index $\maj(P)$ is the sum of $x+y$, where $(x,y)$ ranges over the
right-up corners of the path. The theorem says that 
\begin{equation*}
\Cat(n)_q=\sum_{P\in\Dyck(n)} q^{\maj(P)}.
\end{equation*}
For example, Figure~\ref{fig:major_index} illustrates the expansion
$\Cat(3,4)_q=1+q^2+q^3+q^4+q^6$.

More recently, these $q$-Catalan numbers have been generalized
to the rational case. We saw in Section~\ref{sec:rational_cat}
that the binomial coefficient $\tbinom{a-1+b}{a-1}$ is divisible by
$a$ whenever $\gcd(a,b)=1$, and we called the quotient a rational
Catalan number $\Cat(a,b)=\tbinom{a-1+b}{a-1}/a$. It turns out that the
$q$-binomial coefficient $\tsqbinom{a-1+b}{a-1}_q$ is also divisible by
the $q$-integer $[a]_q$ in the ring $\ZZ[q]$. We learned the following simple proof from Theresia Eisenk\"olbl. First note that 
\begin{equation*}
[a]_q \,\sqbinom{a+b}{a}_q = [a+b]_q \,\sqbinom{a-1+b}{a-1}_q.
\end{equation*}
If $\gcd(a,b)=1$ then we also have $\gcd(a+b,a)=1$. But for any integer $n\geq 1$ recall that $[n]_q$ has irreducible factors given by the cyclotomic polynomials $\Phi_d(q)$ with $d|n$ and $d>1$. It follows that $[a]_q$ and $[a+b]_q$ are coprime, and hence $[a]_q$ divides $\tsqbinom{a-1+b}{a-1}_q$ in the ring $\ZZ[q]$.

Thus we can define the {\em rational $q$-Catalan number}
\begin{equation*}
\Cat(a,b)_q := \frac{1}{[a]_q} \sqbinom{a-1+b}{a-1}_q\in\ZZ[q] \quad \text{when
$\gcd(a,b)=1$.}
\end{equation*}
It turns out that the coefficients are non-negative, though this is much more difficult to prove. As far as we know, the first proof that $\Cat(a,b)_q\in\NN[q]$ was given by Haiman \cite[Propositions 2.5.2--4]{haiman}.
He gave a sketch of a proof showing that positivity would follow from the existence of a certain ``homogeneous system of parameters'' for the polynomial ring $\QQ[x_1,\ldots,x_b]$. In an early draft of the paper he merely conjectured that such an h.s.o.p.\ exists; in the final draft he gave a proof of existence based on an idea of H.~Kraft. Reiner, Stanton and White \cite[Corollary 10.4]{rsw} gave another proof based on the unimodality of the coefficients of $\tsqbinom{a-1+b}{a-1}_q$, which itself is not easy to prove.

Thus it became a problem to find a combinatorial interpretation for the coefficients of $\Cat(a,b)_q$. Since $\Cat(a,b)$ is the number of Dyck paths in an $a\times b$
rectangle, it is natural to search for a statistic
$\Dyck_{a,b}\to\NN$ that restricts to the major index when $b=a+1$
and whose generating function is $\Cat(a,b)_q$. So far no one has been
able to do this.\footnote{The author first learned of this problem
from Dennis Stanton around 2005.}

Instead, our best understanding of the rational $q$-Catalan numbers is currently based on the more subtle $q,t$-Catalan numbers. These
are a family of two-variable polynomials $\Cat(a,b)_{q,t}\in\NN[q,t]$
depending on a pair of coprime integers $\gcd(a,b)=1$, which satisfy
the following properties:\footnote{The specialization
$\Cat(a,b)_{q,1}$ is the $q$-generating function for the area
statistic on rational Dyck paths. This is a different kind of
``rational $q$-Catalan number'' that is much easier to study.}
\begin{itemize}
\item $\Cat(a,b)_{1,1}=\Cat(a,b)$,
\item $\Cat(a,b)_{q,t}=\Cat(b,a)_{t,q}$,
\item $\Cat(a,b)_{q,t}=\Cat(a,b)_{t,q}$,
\item $\Cat(a,b)_q = q^{(a-1)(b-1)/2}\Cat(a,b)_{q,q^{-1}}$.
\end{itemize}
We give a brief history of these polynomials. In 1993, Mark Haiman presented a series of conjectures based on computer experiments  in invariant theory \cite{haiman}. Given a permutation $\sigma \in S_n$ and a polynomial $f(\x,\y)\in\QQ[x_1,\ldots,x_n,y_1,\ldots,y_n]$ in two sets of variables, we define the polynomial $f^\sigma(\x,\y)$ by
\begin{equation*}
f^\sigma(x_1,\ldots,x_n,y_1,\ldots,y_n) := f(x_{\sigma(1)},\ldots,x_{\sigma(n)},y_{\sigma(1)},\ldots,y_{\sigma(n)}).
\end{equation*}
Let $I$ be the ideal generated by non-constant polynomials satisfying $f^\sigma=f$ for all $\sigma\in S_n$, called {\em diagonal invariants}. The quotient ring is called the ring of {\em diagonal coinvariants}:
\begin{equation*}
\DR_n := \QQ[x_1,\ldots,x_n,y_1,\ldots,y_n]/I.
\end{equation*}
This is a finite dimensional $\QQ$-algebra bigraded by $x$-degree and $y$-degree. Let $(\DR_n)_{i,j}$ be the subspace of polynomials with $x$-degree $i$ and $y$-degree $j$. Since the ideal $I$ is bi-homogeneous and $S_n$-invariant, each vector space $(\DR_n)_{i,j}$ carries a representation of $S_n$. If $c_{ij}$ is the multiplicity of the sign character in $(\DR_n)_{i,j}$ then the $q,t$-Catalan number is defined as
\begin{equation*}
\Cat(n)_{q,t} := \sum_{i,j} c_{ij} q^i t^j.
\end{equation*}
The name was inspired by the conjecture that the evaluation at $q=t=1$ yields the classical Catalan number.\footnote{Haiman credits this observation to Stanley.} The paper \cite{haiman} made combinatorial conjectures about the specializations at $t=1$ and $t=1/q$, but did not make a conjecture about the coefficient $c_{ij}$; Haiman said that ``it would be very interesting to have a plausible combinatorial conjecture as to its value''. Since the polynomial $\Cat(n)_{q,1}$ was conjectured to be the generating function for the area statistic of Dyck paths, it was natural to search for another statistic on Dyck paths explaining the $t$-degree. In late 1999, Jim Haglund found such a statistic, which he called ``bounce'', and shared it with Adriano Garsia. After working together for a few months on the conjecture they became pessimistic about proving it. Garsia shared with Haiman that a statistic existed, without giving any details. Within two weeks Haiman also found a statistic, later called ``dinv''.\footnote{Jim Haglund, personal communication.} It turned out that the bounce and dinv statistics are not the same, but Haglund and Haiman found an explicit bijection $\zeta:\Dyck_{n,n}\to\Dyck_{n,n}$ with the properties
\begin{align*}
\dinv(P) &= \area(\zeta(P)),\\
\bounce(P) &= \area(\zeta^{-1}(P)).
\end{align*}
See \cite[Theorem 3.15]{haglund} for details. So, in essence, there is really just the area statistic and a strange bijection $\zeta$. In late 2000, Garsia and Haglund \cite{garsiahaglund} succeeded in proving that
\begin{equation*}
\Cat(n)_{q,t} = \sum_{P\in\Dyck_{n,n}} q^{\area(P)} t^{\bounce(P)}.
\end{equation*}
In particular, since the algebraic definition of $\Cat(n)_{q,t}$ is clearly symmetric in $q$ and $t$, this implies that there must exist a bijection from $\Dyck_{n,n}$ to itself that interchanges the statistics $\area$ and $\bounce$. But Garsia and Haglund's proof is recursive and it does not provide such a bijection; indeed, this is still an open problem. The $\area$ and $\dinv$ statistics were later generalized to ``parking functions'' (which we won't define here), yielding a conjectural combinatorial description of the full bigraded character of $\DR_n$, the so-called {\em shuffle conjecture} \cite{hhlru}. After ten years, the conjecture was proved by Carlsson and Mellit \cite{cm}. As with Garsia and Haglund's proof of the $q,t$-Catalan case, Carlsson and Mellit's proof is recursive and does not ``explain'' the $q,t$-symmetry.

\begin{figure}[h]
\begin{center}
	\includegraphics[scale=1.2]{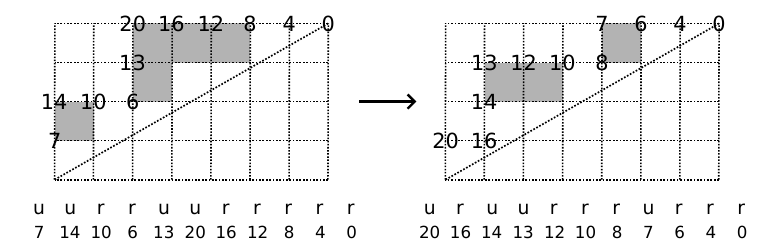}
\end{center}
\caption{The sweep map.}
\label{fig:sweep_map}
\end{figure}

Loehr \cite{loehr} generalized the $\bounce$ and $\dinv$ statistics to Dyck paths in an $n\times mn$ rectangle. For the full rational case, however, it turns out to be easier to generalize the bijection $\zeta$ instead of the statistics. This generalization was discovered at least three independent times: (1) by Loehr and Warrington in terms of Dyck paths \cite{lw}, (2) by Gorsky and Mazin in terms of numerical semigroups \cite{gm}, (3) by Armstrong, Hanusa and Jones in terms of simultaneous $a,b$-core partitions \cite{ahj}. The paper \cite{alw} showed that the three definitions are equivalent. The easiest to describe is the {\em sweep map} formulation. Let $\gcd(a,b)=1$ and consider a Dyck path $P\in\Dyck_{a,b}$. We label each step with an integer as follows: Start with $0$. For each up step add $b$ and for each right step subtract $a$. Finally, define the path $\sweep(P)\in\Dyck_{a,b}$ by sorting the steps by decreasing order on labels. For example, Figure \ref{fig:sweep_map} shows that the path $uurruurrrrr\in\Dyck_{4,7}$ gets sent to\footnote{It is certainly not obvious how to invert the sweep map. This was solved by Thomas and Williams~\cite{tw_sweep}, though the solution is quite complicated.}
\begin{equation*}
\sweep(uurruurrrrr) = uruurrrurrr.
\end{equation*}
 Finally, we use the sweep map to define the {\em rational $q,t$-Catalan numbers}
\begin{equation*}
\Cat(a,b)_{q,t} := \sum_{P\in\Dyck_{a,b}} q^{\area(P)} t^{\area(\sweep(P))},
\end{equation*}
where $\area(P)$ is the number of full squares between the path and the diagonal. Figure \ref{fig:sweep_map} shows that $\area(P)=5$ and $\area(\sweep(P))=3$ for the path $P=uurruurrrrr$.

Each of the papers \cite{ahj,gm,lw} conjectured the symmetry $\Cat(a,b)_{q,t}=\Cat(a,b)_{t,q}$ and the specialization
\begin{align*}
\Cat(a,b)_q &= q^{(a-1)(b-1)/2} \Cat(a,b)_{q,1/q},  \\
\frac{1}{[a]_q} \sqbinom{a-1+b}{a-1}_q &=
\sum_{P\in\Dyck_{a,b}} q^{\area(P)-\area(\sweep(P))+(a-1)(b-1)/2}.
\end{align*}
This time, however, the combinatorial definition of $\Cat(a,b)_{q,t}$ was not preceded by an algebraic definition, so the methods of proof pursued by Garsia, Haglund and Haiman in the classical cases were not immediately available. The eventual proof came when the shuffle conjecture was generalized by Bergeron, Garsia, Leven and Xin \cite{bglx} and then this {\em rational shuffle conjecture} was proved by Mellit \cite{mellit}, generalizing the ideas from Carlsson and Mellit \cite{cm}.

This is still an active area of research, so we end the discussion here. For the purpose of this paper we are concerned with the fact that the statistic  $\Dyck_{a,b}\to\NN$ defined by
\begin{equation*}
P\mapsto \area(P)-\area(\sweep(P))+(a-1)(b-1)/2
\end{equation*}
is difficult to work with and does not shed light on Conjecture \ref{conj:monotone}. In particular, this statistic does not respect the fact that the sets of Dyck paths with fixed $a$ are nested:
\begin{equation*}
\Dyck_{a,0} \subseteq \Dyck_{a,1}\subseteq \Dyck_{a,2} \subseteq
\cdots.
\end{equation*}
We can think of the union $\Dyck_{a,\infty}:= \bigcup_b \Dyck_{a,b}$ as the set of Dyck paths in an $a\times \infty$ rectangle. It would be nice to have a statistic $\Dyck_{a,\infty}\to\NN$ that respects this filtration.

Our idea is to look at lattice points instead of Dyck paths. In this language, we want to find a statistic $\Jsf:\R\to\NN$ on the root lattice that respects the filtration by dilations of the fundamental alcove,
\begin{equation*}
(\R\cap 0\Delta) \subseteq (\R\cap 1\Delta)\subseteq (\R\cap 2\Delta) \subseteq\cdots,
\end{equation*}
and restricts to the rational $q$-Catalan numbers:
\begin{equation*}
\Cat(a,b)_q = \sum_{\x\in (\R\cap b\Delta)} q^{\Jsf(\x)}.
\end{equation*}

\section{Lattice points and $q$-binomial coefficients}
\label{sec:qbinomial}

\begin{figure}[h]
\begin{center}
\includegraphics[scale=0.8]{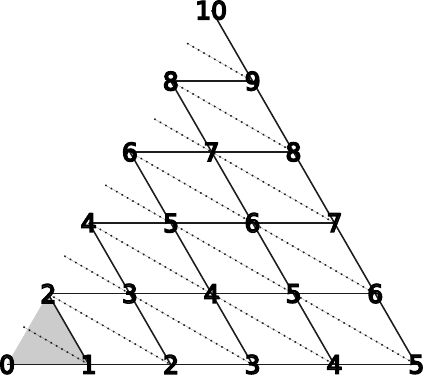}
\end{center}
\caption{Tilted partial order on the set $\L\cap 5\Delta$ when $a=3$.}
\label{fig:qbinomial_lattice}
\end{figure}

In order to interpret $\Cat(a,b)_q$ in terms of lattice points we must
first interpret the $q$-binomial coefficient $\tsqbinom{a-1+b}{a-1}_q$
in terms of lattice points. I claim that this is just
the generating function for the tilted height on the set of lattice
points in the dilated simplex $b\Delta$:
\begin{equation*}
\sum_{\x\in (\L\cap b\Delta)} q^{\Tsf(\x)} = \sqbinom{a-1+b}{a-1}_q.
\end{equation*}
For example, Figure~\ref{fig:qbinomial_lattice} displays the tilted
partial order on the set $\L \cap 5\Delta$ when $a=3$. One can observe
that the tilted height generating function agrees with the expansion
\begin{equation*}
\sqbinom{3-1+5}{3-1}_q=1+1q+2q^2+2q^3+3q^4+3q^5+3q^6+2q^7+2q^8+1q^9+1q^{10}.
\end{equation*}

This fact is not difficult to prove. We give two different proofs
for the sake of exposition.

\begin{figure}[h]
\begin{center}
\includegraphics[scale=0.9]{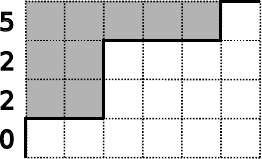}
\end{center}
\caption{A Young diagram.}
\label{fig:young_diagram}
\end{figure}

First we observe that the tilted partial order on $\L\cap b\Delta$ is
isomorphic to Stanley's poset $L(a-1,b)$ of Young diagrams that fit inside
an $(a-1)\times b$ rectangle (see~\cite[Chapter 6]{stanley_book}). A
{\em Young diagram} (also called an {\em integer partition}) is a
weakly-decreasing sequence of integers $y_1\geq y_2\geq \cdots \geq
0$, only finitely many of which are nonzero. The Young diagrams that
fit inside an $(a-1)\times b$ rectangle have the form
\begin{equation*}
L(a-1,b)=\{(y_1,\ldots, y_{a-1})\in\ZZ^{a-1}: b\geq y_1\geq \cdots
\geq y_{a-1}\geq 0\}.
\end{equation*}
We visualize each element of $L(a-1,b)$ as the ``shape'' of the
cells above a lattice path from $(0,0)$ to $(b,a-1)$. For example,
Figure \ref{fig:young_diagram} illustrates the Young diagram $(5,2,2,0)\in
L(4,6)$ as the ``shape'' above a certain lattice path in the 
$(5-1)\times 6$ rectangle.

This bijection between Young diagrams and lattice paths shows that
$\#L(a-1,b)=\tbinom{a-1+b}{a-1}$. The partial order on $L(a-1,b)$
is defined as the componentwise order on the corresponding
integer vectors $\y=(y_1,\ldots,y_{a-1})$, which we can think of
as ``inclusion of diagrams''. Thus the rank of the diagram $\y$
in the poset $L(a-1,b)$ is its {\em area} $y_1+\cdots +y_{a-1}$,
which shows that the rank-generating function of $L(a-1,b)$ is
$\tsqbinom{a-1+b}{a-1}_q$. To prove that $L(a-1,b)$ is isomorphic to
the tilted order on $\L\cap b\Delta$, recall from Section \ref{sec:tilted} that the tilted order
is defined as the componentwise order on coefficients in the tilted
basis. Given a point $\x=x_1\omega_1+\cdots+x_{a-1}\omega_{a-1}$
in standard coordinates, let
$\x=y_1\omega_1+y_2(\omega_2-\omega_1)+\cdots+y_{a-1}(\omega_{a-1}-\omega_{a-2})$
be the expansion in the tilted basis, so that $x_{a-1}=y_{a-1}$ and
$x_i=y_i-y_{i+1}$ for $1\leq i\leq a-2$. In these new coordinates
we have
\begin{align*}
\L\cap b\Delta &= \bigl\{(x_1,\ldots,x_{a-1})\in\ZZ^{a-1}:
\textstyle{\sum_{i=1}^{a-1}} x_i\leq b \text{ and } x_i\geq 0 \text{
for all } 1\leq i\leq a-1\bigr\}\\
&=  \{(y_1,\ldots,y_{a-1})\in\ZZ^{a-1}: y_{a-1}\leq b \text{ and }
y_i-y_{i+1}\geq 0 \text{ for all } 1\leq i\leq a-2\}\\
&=  \{(y_1,\ldots,y_{a-1})\in\ZZ^{a-1}: b\geq y_{a-1}\geq \cdots \geq
y_1\geq 0\},
\end{align*}
which is just the set $L(a-1,b)$ of Young diagrams.

For our second proof we give a generating function argument from Ehrhart theory. Recall that $K$ is the positive cone inside the real vector space spanned by $\L$.  We are interested in the set of positive lattice points $\L\cap K$:
\begin{equation*}
\L\cap K = \{(x_1,\ldots,x_{a-1})\in\ZZ^{a-1}: 0\leq x_i \text{
for all $i$\,}\}.
\end{equation*}
We also consider the higher-dimensional lattice
\begin{equation*}
X = \{(x_0,x_1,\ldots,x_{a-1})\in\ZZ^{a}: 0\leq x_i \text{ for
all $i$\,}\}.
\end{equation*}
We write $\x=(x_1,\ldots,x_{a-1})$ for elements of
$\L\cap K$ and $\hat\x=(x_0,\ldots,x_{a-1})$ for elements 
of~$X$. In Section~\ref{sec:rational_cat} we observed that the projection
$(x_0,x_1,\ldots,x_{a-1})\mapsto (x_1,\ldots,x_{a-1})$ defines a
bijection from $X_b:=\{\hat\x\in X: x_0+\cdots +x_{a-1}=b\}$
to $\L\cap b\Delta=\{\x\in K\cap\L: x_1+\cdots+x_{a-1}\leq
b\}$. Now let $z_0,\ldots,z_{a-1}$ be indeterminates and observe that
we have a geometric series:
\begin{equation*}
\frac{1}{(1-z_0)(1-z_1)\cdots(1-z_{a-1})} = \sum_{\hat\x\in X}
z_0^{x_0}\cdots z_{a-1}^{x_{a-1}}.
\end{equation*}
Substituting $z_i\mapsto tz_i$ for all $i$ gives
\begin{equation*}
\frac{1}{(1-tz_0)(1-tz_1)\cdots(1-tz_{a-1})}= \sum_{\hat\x\in X}
z_0^{x_0}\cdots z_{a-1}^{x_{a-1}}t^{x_0+\cdots + x_{a-1}},
\end{equation*}
and then setting $z_0=1$ gives
\begin{align*}
\frac{1}{(1-t)(1-tz_1)\cdots(1-tz_{a-1})}&=\sum_{\hat\x\in X}
z_1^{x_1}\cdots z_{a-1}^{x_{a-1}}t^{x_0+\cdots +x_{a-1}} \\
&= \sum_{b\geq 0} t^b \sum_{\hat\x\in X_b} z_1^{x_1}\cdots
z_{a-1}^{x_{a-1}}
= \sum_{b\geq 0} t^b \sum_{\x\in (\L\cap b\Delta)} z_1^{x_1}\cdots
z_{a-1}^{x_{a-1}}.
\end{align*}
Finally, substituting $z_i\mapsto q^i$ for all $i$ gives
\begin{equation*}
\frac{1}{(1-t)(1-qt)\cdots(1-q^{a-1}t)} = \sum_{b\geq 0} t^b
\sum_{\x\in (\L\cap b\Delta)} q^{x_1+2x_2+\cdots+(a-1)x_{a-1}}
 = \sum_{b\geq 0} t^b \sum_{\x\in (\L\cap b\Delta)} q^{\Tsf(\x)}.
\end{equation*}
On the other hand, we have the following standard $q$-binomial identity:\footnote{See Stanley \cite{EC1} between (1.86) and (1.87).}
\begin{equation*}
\frac{1}{(1-t)(1-qt)\cdots(1-q^{a-1}t)} = \sum_{b\geq 0} t^b
\sqbinom{a-1+b}{a-1}_q.
\end{equation*}
This completes the second proof.

We showed in Section~\ref{sec:rational_cat} that $\#(R_k\cap
b\Delta)=\Cat(a,b)$ for any $k$ and $\gcd(a,b)=1$. It is certainly
not the case that the tilted height generating function on $\R_k\cap
b\Delta$ is $\Cat(a,b)_q$. For example, if we write $(S)_q$ for
the tilted height generating function on the set $S\subseteq\L$
then from Figure~\ref{fig:qbinomial_lattice} we can see that
\begin{align*}
(\R_0\cap 5\Delta)_q &=
1+2q^3+3q^6+q^9,\\
(\R_1\cap 5\Delta)_q &=
q+3q^4+2q^7+q^{10},\\
(\R_2\cap 5\Delta)_q &=
2q^2+3q^5+2q^8,
\end{align*}
each of which is a $q$-analogue of the number $\Cat(3,5)=7$. But the
$q$-Catalan number is
\begin{equation*}
\Cat(3,5)_q = 1+q^2+q^3+q^4+q^5+q^6+q^8.
\end{equation*}
In the next section we attempt to modify the tilted height to obtain a statistic whose generating function on the set $\R\cap b\Delta$ is $\Cat(a,b)_q$.

\section{Johnson statistics}
\label{sec:johnson}

At the end of Section~\ref{sec:rational_qcat} we suggested that there should be a statistic $\Jsf:\R\to\ZZ$ on the
full root lattice $\R$ with the property that $\Cat(a,b)_q=\sum_{\x\in
(\R\cap b\Delta)} q^{\Jsf(\x)}$ for all integers $b\geq 1$ satisfying
$\gcd(a,b)=1$. The existence of such a statistic is closely related to a conjecture of
Paul Johnson, which we now describe. Johnson mentioned \cite[Section 4.4 and Remark 4.8]{johnson} that this conjecture
was motivated partly by Chen--Ruan cohomology, and partly by the
following rearrangement of the rational $q$-Catalan number:\footnote{Johnson's statement of this formula has a typo.}
\begin{equation*}
\Cat(a,b)_q
=\frac{(1-q^{b+1})(1-q^{b+2})\cdots(1-q^{b+a-1})}{(1-q^a)(1-q^{2a})\cdots (1-q^{a(a-1)})}
[a]_{q^2}[a]_{q^3}\cdots [a]_{q^{a-1}}.
\end{equation*}
Recall that we write $\L=\ZZ^{a-1}$ for the weight lattice and  $\R\subseteq\L$ for the root lattice whose points satisfy $x_1+2x_2+\cdots+(a-1)x_{a-1} \equiv 0$ mod $a$. Note that the sublattice $a\L$ is contained in $\R$. The following observation is crucial for this section:
\begin{equation*}
\#(\R/a\L) = a^{a-2}.
\end{equation*}
Indeed, since $\L$ has rank $a-1$ we must have $\#(\L/c\L)=c^{a-1}$ for any integer $c\geq 1$. In particular, we have $\#(\L/a\L)=a^{a-1}$. Furthermore, we recall from Section \ref{sec:lattices} that $\L/\R\cong \ZZ/a\ZZ$, and hence $\#(\L/\R)=a$. Then the tower law gives $\#(\L/\R)\cdot \#(\R/a\L)=\#(\L/a\L)$.

The following conjecture appears in the arxiv version of \cite{johnson}, but not in the published version.

\noindent{\bf Johnson's conjecture.} \cite[Conjecture
4.6]{johnson}. There exists a statistic $\iota:\R/a\L\to\NN$ on the
set of cosets of $a\L$ in $\R$ satisfying the following two properties:
\begin{itemize}
\item The generating function for $\iota$ is
\begin{equation*}
\sum_{\mathfrak{c}\,\in (\R/a\L)} q^{\iota(\mathfrak{c})} =
[a]_{q^2}[a]_{q^3}\cdots [a]_{q^{a-1}}.
\end{equation*}
\item When $\gcd(a,b)=1$ we have
\begin{equation*}
\Cat(a,b)_q = \sum_{\mathfrak{c}\,\in(\R/a\L)} q^{\iota(\mathfrak{c})}
\sqbinom{a-1+b/a-s(\mathfrak{c},b)}{a-1}_{q^a}
\end{equation*}
for some small rational numbers $s(\mathfrak{c},b)\in\QQ$ depending
on $\mathfrak{c}$ and $b$.
\end{itemize}
Johnson suggested that $\iota$ should be an ``age'' function for the Chen--Ruan cohomology of an orbifold, but he admitted that ``this discussion is rather vague''. We will give a more concrete combinatorial interpretation.

Recall the definitions of the fundamental (half-open) parallelotope and fundamental box:
\begin{align*}
\Pi &= \{(x_1,\ldots,x_{a-1})\in\RR^{a-1}: 0\leq x_i< 1 \text{
for all $i$\,}\} \subseteq \L\otimes\RR,\\
\Box = a\Pi &= \{(x_1,\ldots,x_{a-1})\in\RR^{a-1}: 0\leq x_i< a \text{
for all $i$\,}\} \subseteq \L\otimes\RR.
\end{align*}
For any integer $b\geq 0$ we note that $\#(\L\cap b\Pi)=b^{a-1}$. In particular, we have $\#(\L\cap \Box)=a^{a-1}$. We will show that $\L\cap\Box$ is a set of coset representatives for $\L/\a\L$ and $\R\cap \Box$ is a set of coset representatives for $\R/a\L$. To see this, note that each lattice point $\x\in\L$ has a unique
quotient $(\x\text{ quo } a)\in\L$ and remainder $(\x \text{ rem
}a)\in\L$ such that $\x=a(\x\text{ quo } a)+(\x \text{ rem } a)$
and the coordinates of $(\x \text{ rem } a)=(r_1,\ldots,r_{a-1})$
satisfy $0\leq r_i\leq a-1$ for all $i$. Furthermore, we observe
that $\x\in\R$ if and only if $(\x \text{ rem } a)\in\R$. Indeed, for any $\x\in\L$ we have $\x\in\R$ if and only if $\Tsf(\x)\equiv 0$ mod $a$. But since the tilted height is linear we have
\begin{align*}
\Tsf(\x) &= \Tsf[a(\x\text{ quo } a)+(\x \text{ rem } a)]\\
&=a\Tsf(\x\text{ quo } a)+\Tsf(\x \text{ rem } a) \\
&\equiv \Tsf(\x \text{ rem } a) \text{ (mod $a$)}.
\end{align*}
Thus the function $\x\mapsto (\x \text{ rem } a)$ defines two bijections:
\begin{align*}
\L/a\L &\quad\to \quad \L\cap \Box,\\
\R/a\L &\quad\to \quad \R\cap \Box.
\end{align*}
In particular, we observe that $\#(\R\cap \Box)=\#(\R/a\L)=a^{a-2}$. Our key idea is to replace Johnson's ``age'' statistic $\R/a\L\to\ZZ$ by a statistic on the root lattice $\R\to\ZZ$ with a certain periodicity mod $a$.

\begin{definition}\label{def:johnson_statistic} Consider a function
$\Jsf:\R\to\ZZ$ on the root lattice of type $A_{a-1}$. We say 
that~$\Jsf$ is a {\em Johnson statistic} if it satisfies the following two properties:
\begin{itemize}
\item (Periodic) For any $\x\in \R$ and $\y\in\L$ we have
$\Jsf(\x+a\y)=\Jsf(\x)+a\Tsf(\y)$. That is,
\begin{equation*}
\Jsf(x_1+ay_1,\ldots,x_{a-1}+ay_{a-1}) = \Jsf(x_1,\ldots,x_{a-1})+a(y_1+2y_2+\cdots+(a-1)y_{a-1}).
\end{equation*}
\item (Catalan) For any integer $b\geq 0$ coprime to $a$ we have
\begin{equation*}
\Cat(a,b)_q = \sum_{\x \in (\R\cap b\Delta)} q^{\Jsf(\x)}.
\end{equation*}
\end{itemize}

\end{definition}

\begin{conjecture}\label{conj:johnson_statistics} There exists at least one Johnson statistic for each value of $a$.
\end{conjecture}

\begin{figure}[h]
\begin{center}
\includegraphics[scale=0.9]{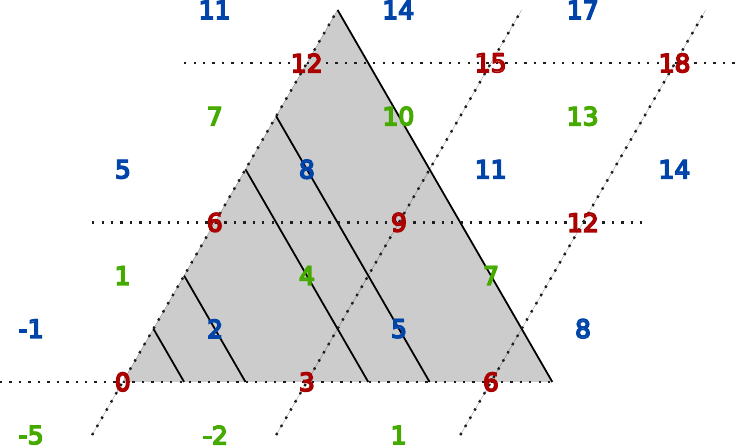}
\end{center}
\caption{The unique Johnson statistic when $a=3$.}
\label{fig:johnson_a2}
\end{figure}

By periodicity, a Johnson statistic is determined by its values on the coset representatives $\R\cap\Box$. In the case $a=3$ this is the set
\begin{align*}
\R\cap\Box &= \{(x_1,x_2)\in\ZZ^2: 0\leq x_1,x_2 <3 \text{ and } x_1+2x_2 \equiv 0 \text{ mod $3$}\}\\
 &= \{(0,0), (1,1), (2,2)\}.
\end{align*}
In this case we will see in Section \ref{sec:ribbons} that there exists a unique Johnson statistic, which is defined by $\Jsf(0,0)=0$, $\Jsf(1,1)=2$ and $\Jsf(2,2)=4$. Then we extend this by periodicity:
\begin{equation*}
\Jsf(x_1+3y_1,x_2+3y_2)=\Jsf(x_1,x_2)+3(y_1+2y_2).
\end{equation*}
Figure \ref{fig:johnson_a2} shows this statistic. The three colors label the $3=3^{3-2}$ cosets in $\R/3\L$.

Translations of the fundamental box are shown by dotted lines. The solid diagonal lines show dilations of the fundamental alcove $b\Delta$ for various integers $b$ with $\gcd(3,b)=1$. Observe that the sum of $q^{\Jsf(\x)}$ over points $\x \in (\R\cap b\Delta)$ gives the rational $q$-Catalan numbers:
\begin{align*}
\Cat(3,1)_q&=1,\\
\Cat(3,2)_q&=1+q^2,\\
\Cat(3,4)_q&=1+q^2+q^3+q^4+q^6,\\
\Cat(3,5)_q&=1+q^2+q^3+q^4+q^5+q^6+q^8,\\
\Cat(3,7)_q&=1+q^2+q^3+q^4+q^5+2q^6+q^7+q^8+q^9+q^{10}+q^{12}.
\end{align*}
Note that this explains monotonicity. That is, Conjecture \ref{conj:johnson_statistics} implies Conjecture \ref{conj:monotone}.

Next we show that Conjecture \ref{conj:johnson_statistics} implies Johnson's conjecture. We  separate this into two parts, each of which is based on a certain geometric partition. The first part is a bit easier.

\begin{theorem}\label{thm:johnson1}
Let $\Jsf:\R\to\ZZ$ be a Johnson statistic and let $\gcd(a,b)=1$. Then we have
\begin{equation*}
\Cat(a,b)_q = \sum_{\x\in (\R\cap \Box)} q^{\Jsf(\x)} \sqbinom{a-1+\lfloor
(b-\sum_i x_i)/a\rfloor}{a-1}_{q^a},
\end{equation*}
where the sum is taken over points $\x=(x_1,..,x_{a-1})$ in $\R\cap\Box$. 
\end{theorem}

\begin{proof}
Since $\R$ is the disjoint union of cosets $\x+a\L$ we note that $\R\cap b\Delta$ is the disjoint union of the sets $(\x+a\L)\cap b\Delta$, where $\x$ runs over the points of $\R\cap \Box$. If $\x$ is in $\R\cap\Box$, so that $0\leq x_i\leq a-1$, then we also note that
\begin{align*}
(\x+a\L)\cap b\Delta &= \{\x+a\y:  x_i+ay_i\geq 0 \text{ and }  \sum_i (x_i+ay_i)\leq b\} \\
&= \{\x+a\y: y_i\geq 0 \text{ and } \sum_i y_i \leq \lfloor (b-\sum_i x_i)/a\rfloor \} \\
&= \x+a(\L\cap b'\Delta),
\end{align*}
where we write $b'=\lfloor (b-\sum_i x_i) /a\rfloor$ to save space. By the periodicity and Catalan properties of the Johnson statistic $\Jsf$ we have
\begin{align*}
\sum_{\z\in ((\x+a\L)\cap b\Delta)} q^{\Jsf(\z)} &= \sum_{\z\in (\x+a(\L\cap b'\Delta))} q^{\Jsf(\z)}\\
&= \sum_{\y\in (\L\cap b'\Delta)} q^{\Jsf(\x+a\y)}\\
&= \sum_{\y\in (\L\cap b'\Delta)} q^{\Jsf(\x)+a\Tsf(\y)}\\
&= q^{\Jsf(\x)} \sum_{\y\in (\L\cap b'\Delta)} (q^a)^{\Tsf(\y)}\\
&= q^{\Jsf(\x)} \sqbinom{a-1+b'}{a-1}_{q^a}.
\end{align*}
Then summing over cosets gives
\begin{equation*}
\Cat(a,b)_q =\sum_{\x\in(\R\cap b\Delta)} q^{\Jsf(\x)}= \sum_{\substack{ \x\in (\R\cap\Box) \\ \z\in ((\x+a\L)\cap b\Delta)}} q^{\Jsf(\z)}= \sum_{\x\in (\R\cap\Box)} q^{\Jsf(\x)} \sqbinom{a-1+b'}{a-1}_{q^a}.
\end{equation*}
\end{proof}

This coset decomposition is illustrated by the colors in Figure \ref{fig:johnson_a2}. In this example the sum over the three coset representatives gives
\begin{equation*}
\Cat(3,b)_q = \sqbinom{2+\lfloor (b-0)/3\rfloor}{2}_{q^3}+q^2\sqbinom{2+\lfloor
(b-2)/3\rfloor}{2}_{q^3}+q^4\sqbinom{2+\lfloor
(b-4)/3\rfloor}{2}_{q^3}
\end{equation*}
for any $\gcd(3,b)=1$, which agrees with Example 4.10 in Johnson \cite{johnson}.

For the other part of Johnson's conjecture we need a lemma. Christian Krattenthaler pointed out to us that this identity follows from the $q$-binomial theorem, but it doesn't seem to be a standard result. We will prove a stronger version than we need for Johnson's conjecture since we will use it below to show that Johnson statistics satisfy an analogue of Brion's theorem (see Theorem \ref{thm:brion}).

\begin{lemma}\label{lem:johnson} Define the $q$-Pochhammer symbol $(u;q)_n = (1-u)(1-uq)\cdots (1-uq^{n-1})$. For any set $J\subseteq\NN$ we let $\#J$ denote the cardinality and let $\sum J$ denote the sum $\sum_{j\in J} j$. Then
\begin{equation*}
\sum_{J\subseteq\{1,\ldots,a-1\}} (-1)^{\#J}q^{a\sum J} \frac{(uq^{-a\#J};q)_{a-1}} {(q;q)_{a-1}} = \prod_{k=1}^{a-1} [a]_{q^k},
\end{equation*}
Note that the result is independent of $u$. More generally, for any fixed $0\leq i\leq a-1$ we  have
\begin{equation*}
\sum_{J\subseteq\{0,1,\ldots,a-1\}\setminus\{i\}} (-1)^{\#J}q^{a\sum J} \frac{(uq^{-a\#J};q)_{a-1}} {(q;q)_{a-1}} = u^i q^{-(a-1)\binom{i+1}{2}-i} \prod_{k=1}^i [a]_{q^k} \prod_{k=1}^{a-1-i} [a]_{q^k}.
\end{equation*}
\end{lemma}

\begin{proof} The $q$-binomial theorem in Pochhammer notation says that
\begin{equation*}
(-x;q)_{a-1} = \sum_{r=0}^{a-1} q^{\binom{r}{2}} \frac{(q;q)_{a-1}}{(q;q)_r (q;q)_{a-1-r}}\,  x^r.
\end{equation*}
Applying this with $x=-uq^{-a\#J}$ gives
\begin{align*}
\sum_{J\subseteq\{0,1,\ldots,a-1\}\setminus\{i\}} (-1)^{\#J}q^{a\sum J} \frac{(uq^{-a\#J};q)_{a-1}} {(q;q)_{a-1}} &= \sum_J (-1)^{\#J}q^{a\sum J}\sum_{r=0}^{a-1} \frac{(-1)^r u^r q^{-ar\#J} q^{\binom{r}{2}} }{(q;q)_r (q;q)_{a-1-r}} \\
&= \sum_{r=0}^{a-1} \frac{(-1)^r u^r q^{\binom{r}{2}} }{(q;q)_r (q;q)_{a-1-r}} \sum_{J} (-1)^{\#J} q^{a\sum J}q^{-ar\#J} \\
&= \sum_{r=0}^{a-1} \frac{(-1)^r u^r q^{\binom{r}{2}} }{(q;q)_r (q;q)_{a-1-r}}\, S_r,
\end{align*}
where
\begin{equation*}
S_r = \sum_{J} (-1)^{\#J} q^{a\sum J}q^{-ar\#J} 
= \prod_{j\in \{0,1,\ldots,a-1\}\setminus\{i\}} (1-q^{a(j-r)}).
\end{equation*}
Note that $S_r=0$ unless $r=i$, in which case
\begin{align*}
S_i &=  \prod_{j\in \{0,1,\ldots,a-1\}\setminus\{i\}} (1-q^{a(j-i)}) \\
&= \prod_{k=1}^i (1-q^{-ak}) \prod_{k=1}^{a-1-i} (1-q^{ak}) \\
&= \prod_{k=1}^i -q^{-ak}(1-q^{ak}) \prod_{k=1}^{a-1-i} (1-q^{ak}) \\
&= (-1)^i q^{-a\binom{i+1}{2}}\prod_{k=1}^i (1-q^{ak}) \prod_{k=1}^{a-1-i} (1-q^{ak}) \\
&= (-1)^i q^{-a\binom{i+1}{2}} (q^a;q^a)_i (q^a;q^a)_{a-1-i}.
\end{align*}
Hence
\begin{align*}
\sum_{J\subseteq\{0,1,\ldots,a-1\}\setminus\{i\}} (-1)^{\#J}q^{a\sum J} \frac{(uq^{-a\#J};q)_{a-1}} {(q;q)_{a-1}} &=  \frac{(-1)^i u^i q^{\binom{i}{2}} }{(q;q)_i (q;q)_{a-1-i}} \, S_i \\
&=  \frac{(-1)^i u^i q^{\binom{i}{2}} (-1)^i q^{-a\binom{i+1}{2}} (q^a;q^a)_i (q^a;q^a)_{a-1}}{(q;q)_i (q;q)_{a-1-i}} \\
&= u^i q^{-(a-1)\binom{i+1}{2}-i}\frac{(q^a;q^a)_i (q^a;q^a)_{a-1}}{(q;q)_i (q;q)_{a-1-i}}\\
&= u^i q^{-(a-1)\binom{i+1}{2}-i}\prod_{k=1}^i [a]_{q^k} \prod_{k=1}^{a-1-i} [a]_{q^k}.
\end{align*}
\end{proof}

\begin{theorem}\label{thm:johnson2}
Let $\Jsf:\R\to\ZZ$ be a Johnson statistic. Then we have
\begin{equation*}
\sum_{\x\in (\R\cap \Box)} q^{\Jsf(\x)} = [a]_{q^2} [a]_{q^3} \cdots [a]_{q^{a-1}}.
\end{equation*}
More generally, suppose that $\gcd(a,b)=1$ and fix an index $1\leq i\leq a-1$. Consider the isometry $\phi_b$ from Section \ref{sec:rational_cat} and apply the $i$-th power of $\phi_b$ to $\Box$ to get the ``rotated''\footnote{The map $\phi_b$ is an isometry of order $a$ but it is not a rotation. It is actually related to the standard Coxeter element in the symmetric group $S_a$ so it becomes a rotation when restricted to the two-dimensional Coxeter plane passing through the centroid of $b\Delta$.} (half-open) parallelotope $\Box_b^i:=\phi_b^i(\Box)$, which satisfies
\begin{equation*}
\Box_b^i =\{(x_1,\ldots,x_{a-1})\in\RR^{a-1}: 0\leq x_j <a \text{ for $j\neq i$ and } b-a< \sum_{j=1}^{a-1}x_j \leq b\}.
\end{equation*}
Then we have
\begin{equation*}
\sum_{\x\in (\R\cap \Box_b^i)} q^{\Jsf(\x)} = q^{bi-(a-1)\binom{i+1}{2}} \frac{1}{[a]_q} \prod_{k=1}^i [a]_{q^k} \prod_{k=1}^{a-1-i} [a]_{q^k}.
\end{equation*}
\end{theorem}

\begin{proof}
Let $\Jsf:\R\to\ZZ$ be a Johnson statistic. For any finite\footnote{In Theorem \ref{thm:brion} below we will extend this to certain infinite sets $S\subseteq\R$.} set $S\subseteq\R$ we will use the shorthand notation
\begin{equation*}
(S)_\Jsf := \sum_{\x\in S} q^{\Jsf(\x)} \in\ZZ[q^\pm].
\end{equation*}
For any two finite subsets $S,T\subseteq\R$ we note that $(S\cup T)_\Jsf = (S)_\Jsf + (T)_\Jsf - (S\cap T)_\Jsf$. We will use this property in an inclusion-exclusion argument.

Let $\gcd(a,b)=1$ and assume that $b\geq a(a-1)$ so that $\Box$ is completely contained in the dilated simplex $b\Delta$. (At the end of the proof we will show how to remove this assumption.) In this case we will express $(\R\cap\Box_b^i)_\Jsf$ as an alternating sum,
\begin{equation*}
(\R\cap\Box_b^i)_\Jsf=\sum_{I\subseteq \{0,\ldots,a-1\}\setminus\{i\}} (-1)^{\#I} (\R\cap\Delta_I)_\Jsf,
\end{equation*}
where each $\Delta_I$ is a simplex of the form $a\y_I + c_I \Delta$ for some $\y_I\in\L$ and $\gcd(a,c_I)=1$.\footnote{The point $\y_I$ and integer $c_I$ also depend on the choice of $b$ and $i$, but we suppress this in the notation to improve readability.} Hence from the periodic and Catalan properties of $\Jsf$ we can write
\begin{equation*}
(\R\cap\Delta_I)_\Jsf = q^{a\Tsf(\y_I)} \Cat(a,c_I)_q.
\end{equation*}
To see how the construction works we first discuss the case $i=0$, where $\Box_b^i = \Box$ is the ``unrotated'' fundamental box. In this case, for each subset $I\subseteq \{1,\ldots,a-1\}$ we let $\y_I$ be the indicator vector whose $i$-th coordinate is $1$ when $i\in I$ and $0$ when $i\not\in I$, so that $\Tsf(\y_I) = \sum I$. Note that the points $a\y_I$ are just the vertices of $\Box$. Next define $\Delta_I = (a\y_I+K)\cap b\Delta$ where $K$ is the positive cone. It is not difficult to check that $\Delta_I=a\y_I + (b-a\#I)\Delta$, so we can take $c_I=b-a\#I$, which satisfies $\gcd(a,c_I)=1$. One can also check that $\Delta_I=\bigcap_{i\in I} \Delta_{\{i\}}$ and that $\R\cap\Box$ is the complement of the union of the sets $\R\cap\Delta_I$ over non-empty $I\subseteq\{1,\ldots,a-1\}$.\footnote{Here is where we use the fact that $b\geq a(a-1)$.} Hence by inclusion-exclusion we obtain
\begin{equation*}
(\R\cap\Box)_\Jsf = \sum_{I\subseteq\{1,\ldots,a-1\}} (-1)^{\#I} q^{a\sum I} \Cat(a,b-a\#I)_q.
\end{equation*}
To complete the proof of the case $i=0$, we express the rational $q$-Catalan number in terms of $q$-Pochhammer symbols as follows:
\begin{equation*}
\Cat(a,b)_q = \frac{1}{[a]_q} \frac{(q^{b+1};q)_{a-1}}{(q;q)_{a-1}}.
\end{equation*}
Then setting $u=q^{b+1}$ in Lemma \ref{lem:johnson} gives 
\begin{align*}
(\R\cap\Box)_\Jsf &= \sum_{I\subseteq\{1,\ldots,a-1\}} (-1)^{\#I} q^{a\sum I} \Cat(a,b-a\#I)_q \\
&=  \frac{1}{[a]_q} \sum_{I\subseteq \{1,\ldots,a-1\}} (-1)^{\#I} q^{a\sum I}  \frac{(q^{b+1-a\#I};q)_{a-1}}{(q;q)_{a-1}}\\
&= \frac{1}{[a]_q} \prod_{k=1}^{a-1} [a]_{q^k} \\
&= [a]_{q^2} \cdots [a]_{q^{a-1}}.
\end{align*}
Thus we have proved the result for $\gcd(a,b)=1$ and $i=0$, which is the case of Johnson's conjecture. The general case  $1\leq i\leq a-1$ is obtained by applying the $i$-th power of the isometry $\phi_b$ to the simplices $\Delta_I$ for $I\subseteq\{1,\ldots,a-1\}$. Let $\tau_i$ be the bijection from subsets of $\{1,\ldots,a-1\}$ to subsets of $\{0,\ldots,a-1\}\setminus\{i\}$ defined on elements by $\tau_i(j)=i+j$ mod $a$. If $J=\tau_i(I)$ then one can check that the lowest point of $\phi_b^i(\Delta_I)$ is $a\y_{J\setminus \{0\}}$, with height $T(a\y_{J\setminus\{0\}})=a\sum J$. Since $\phi_b$ is an isometry, the dilation factors of the simplices $\phi_b^i(\Delta_I)$ and $\Delta_I$ are the same; namely, $b-a\#I=b-a\#J$. In other words, we have $\phi_b^i(\Delta_I) = a\y_{J\setminus\{0\}} + (b-a\#J)\Delta$. Since summing over $I\subseteq \{1,\ldots,a-1\}$ is the same as summing over $J=\tau_i(I)\subseteq\{0,\ldots,a-1\}\setminus\{i\}$ we get
\begin{equation*}
(\R\cap\Box_b^i)_\Jsf= \sum_{J\subseteq\{0,\ldots,a-1\}\setminus\{i\}} (-1)^{\#J} q^{a\sum J} \Cat(a,b-a\#J)_q.
\end{equation*}
Finally, setting $u=q^{b+1}$ in Lemma \ref{lem:johnson} gives the result. Note that we have proved this result for $\gcd(a,b)=1$ and $b\geq a(a-1)$. To prove the result for all $\gcd(a,b)=1$ and $b\in\ZZ$ (including negative $b$), suppose that $b<a(a-1)$ and choose $k\in\ZZ$ such that $b':=b+ka\geq a(a-1)$. Since $\gcd(a,b')=1$ we know that the result holds for $(\R\cap \Box_{b'}^i)_\Jsf$. But we also observe that $\Box_b^i=\Box_{b'}^i-ka\omega_i$; hence from the periodicity property of $\Jsf$ we have
\begin{align*}
(\R\cap\Box_b^i)_\Jsf &= \sum_{\x\in (\R\cap\Box_{b}^i)} q^{\Jsf(\x)}\\
&= \sum_{\x\in (\R\cap\Box_{b'}^i)} q^{\Jsf(\x-ka\omega_i)}\\
&= \sum_{\x\in (\R\cap\Box_{b'}^i)} q^{\Jsf(\x)}q^{-kai}\\
&=q^{-kai} (\R\cap\Box_{b'}^i)_\Jsf,
\end{align*}
and the result follows from this.
\end{proof}

\begin{figure}[h]
\begin{center}
\includegraphics[scale=0.9]{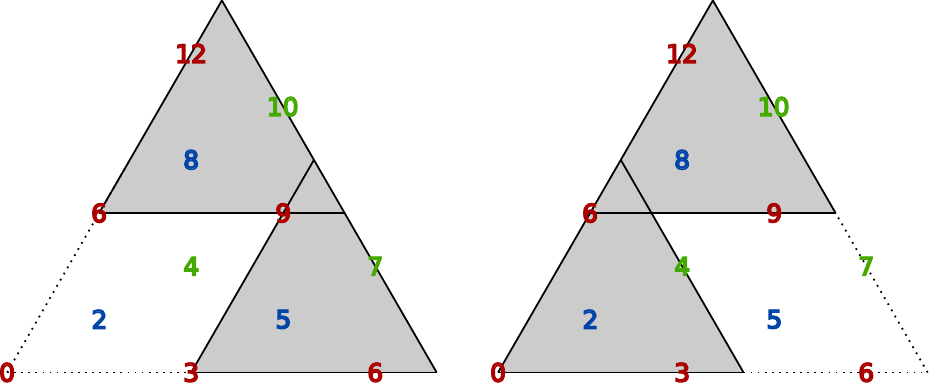}
\end{center}
\caption{Two vertex contributions to $\Cat(3,7)_q$.}
\label{fig:johnson_a2_brion}
\end{figure}

For example, Figure \ref{fig:johnson_a2_brion} demonstrates the inclusion-exclusion construction for the case $(a,b)=(3,7)$ with $i=0$ (left) and $i=1$ (right). To be precise, the picture on the left shows that
\begin{align*}
(\R\cap\Box)_\Jsf &= \sum_{J\subseteq \{1,2\}} (-1)^{\#J} q^{3\sum J} \Cat(3,7-3\# J)_q \\
&= \Cat(3,7)_q - q^3 \Cat(3,4)_q - q^6\Cat(3,4)_q + q^9 \Cat(3,1)_q \\
&= 1+q^2+q^4\\
&= [3]_{q^2},
\end{align*}
and the picture on the right shows that
\begin{align*}
(\R\cap\Box_7^1)_\Jsf &= \sum_{J\subseteq \{0,1,2\}\setminus\{1\}} (-1)^{\#J} q^{3\sum J} \Cat(3,7-3\# J)_q \\
&= \Cat(3,7)_q -  \Cat(3,4)_q - q^6\Cat(3,4)_q + q^6 \Cat(3,1)_q \\
&= q^5+q^6+q^7 \\
&= q^5 [3]_{q}.
\end{align*}

We can use this extra level of generality to show that Johnson statistics satisfy a certain $q$-analogue of Brion's theorem. For brevity we omit some computational details. First let $K=\{(x_1,\ldots,x_{a-1}): 0\leq x_i \text{ for all $i$}\,\}$ be the positive cone and consider the ``rotated'' cones
\begin{equation*}
K_b^i := \phi_b^i(K)= \{ (x_1,\ldots,x_{a-1}): 0\leq x_j \text{ for $j\neq i$ and } \sum_{j=1}^{a-1} x_j\leq b\}.
\end{equation*}
Note that $K_b^i$ is the tangent cone of the simplex $b\Delta$ at the vertex $b\omega_i$, and we have $K_b^i=b\omega_i+\RR_{\geq 0}\{\omega_1^i,\ldots,\omega_{a-1}^i\}$ for the ``rotated'' basis defined by  $\omega_j^i:=\omega_j-\omega_i$ for $j\neq i$ and $\omega_i^i:=-\omega_i$.\footnote{This notation also works for $i=0$ if we define $\omega_0$ to be the zero vector.} The integer points are $\L\cap K_b^i = b\omega_i+\NN\{\omega_1^i,\ldots,\omega_{a-1}^i\}$ and the tilted height generating function for these integer points is a formal Laurent series:
\begin{equation*}
\sum_{\x \in (\L\cap K_b^i)} q^{\Tsf(\x)} = \sum_{\y\in\NN\{\omega_1^i,\ldots,\omega_{a-1}^i\}} q^{\Tsf(b\omega_i +\y)} =\frac{q^{\Tsf(b\omega_i)}}{ \prod_{j=1}^{a-1}(1-q^{\Tsf(\omega_j^i)})} = \frac{q^{bi}}{\displaystyle{\prod_{j\in\{0,\ldots,a-1\}\setminus\{i\}}} (1-q^{j-i})}.
\end{equation*}
The points of the root lattice $\R\cap K_b^i$ can be partitioned by cosets of $a\L$ to get
\begin{equation*}
\R\cap K_b^i = \coprod_{\x\in(\R\cap \Box_b^i)} (\x+\NN\{a\omega_1^i,\ldots,a\omega_{a-1}^i\}),
\end{equation*}
which implies that the Johnson generating function for the $i$-th vertex cone is
\begin{align*}
(\R\cap K_b^i)_\Jsf &= \sum_{\x \in (\R\cap K_b^i)} q^{\Jsf(\x)} \\
&= \sum_{\substack{\x\in(\R\cap\Box_b^i) \\ \y\in\NN\{\omega_1^i,\ldots,\omega_{a-1}^i\}}}q^{\Jsf(\x+a\y)} \\\
&= \sum_{\x\in(\R\cap\Box_b^i)} q^{\Jsf(\x)} \sum_{\y\in\NN\{\omega_1^i,\ldots,\omega_{a-1}^i\}} (q^a)^{\Tsf(\y)}\\
&= \frac{(\R\cap\Box_b^i)_\Jsf}{\displaystyle{\prod_{j\in\{0,\ldots,a-1\}\setminus\{i\}}} (1-q^{a(j-i)})}\\
&= (-1)^i q^{a\binom{i+1}{2}} \frac{(\R\cap\Box_b^i)_\Jsf}{(q^a;q^a)_i(q^a;q^a)_{a-1-i}}.
\end{align*}
Substituting the formula for $(\R\cap\Box_b^i)_\Jsf$ from Theorem \ref{thm:johnson2} and simplifying gives
\begin{equation*}
(\R\cap K_b^i)_\Jsf = (-1)^i q^{bi+\binom{i+1}{2}} \frac{1}{[a]_q} \frac{1}{(q;q)_i(q;q)_{a-1-i}}.
\end{equation*}
Finally, by summing over $i$ and using the $q$-binomial theorem we find that the sum of the formal Laurent series for the vertex cones of $b\Delta$ equals the rational $q$-Catalan number, which is the finite Johnson generating function over $b\Delta$. Thus we obtain an interesting $q$-analogue of Brion's theorem for simplices.

\begin{theorem}\label{thm:brion} Let $\Jsf:\R\to\ZZ$ be a Johnson statistic and let $K_b^i$ be the tangent cone to the dilated fundamental alcove $b\Delta$ at the vertex $b\omega_i$, with $\omega_0:=0$ being the origin. Then we have
\begin{align*}
\sum_{i=0}^{a-1} (\R\cap K_b^i)_\Jsf = \Cat(a,b)_q = (\R\cap b\Delta)_\Jsf,
\end{align*}
where $(S)_\Jsf =\sum_{\x\in S} q^{\Jsf(\x)}$ is the Johnson point enumerator for subsets $S\subseteq\R$.\footnote{It follows from the proof that the formal Laurent series $(\R\cap K_b^i)_\Jsf \in\ZZ[[q^{\pm}]]$ are well defined.}
\end{theorem}

Consider the example $(a,b)=(3,7)$. Summing over the three vertex contributions (two of which are shown in Figure \ref{fig:johnson_a2_brion}) gives 
\begin{equation*}
\frac{1+q^2+q^4}{(1-q^3)(1-q^6)}+\frac{q^5+q^6+q^7}{(1-q^{-3})(1-q^3)}+\frac{q^8+q^{10}+q^{12}}{(1-q^{-6})(1-q^{-3})} = \Cat(3,7)_q.
\end{equation*}
We consider Theorem \ref{thm:brion} to be indirect evidence in favor of Conjecture \ref{conj:johnson_statistics}. More generally, we believe that Theorem \ref{thm:brion} holds for any alcoved polytope (a polytope that is a union of alcoves for the affine Weyl group), but we leave this for future investigation.\footnote{This is not an example of Chapoton's $q$-Ehrhart theory \cite{chapoton} because the statistic $\Jsf$ is not linear.}

\section{$q$-Catalan germs}

When attempting to prove a conjecture it is often helpful to find a stronger conjecture that is still true, thus narrowing the search for a proof. In this section we discuss a refinement of rational $q$-Catalan numbers into more basic polynomials, which we call $q$-Catalan germs, and we conjecture that these more basic polynomials have non-negative coefficients. Then in Section \ref{sec:greedy} we show that this positivity conjecture implies the existence of Johnson statistics.

For any real numbers $\alpha\leq \beta$ we consider a slice of the fundamental box:
\begin{equation*}
\Box[\alpha,\beta]:=\{(x_1,\ldots,x_{a-1})\in \Box: \alpha \leq \sum_i x_i\leq \beta\}.
\end{equation*}
Note that this set is empty when $\beta<0$ or $a(a-1)\leq \alpha$.

\begin{theorem}\label{thm:germs}
Let $m\leq n$ be integers with $\gcd(a,n)=\gcd(a,m-1)=1$. Then the tilted height generating function for the integer points $\L=\ZZ^{a-1}$ in $\Box[m,n]$ is divisible by $[a]_q$ in the ring of integer polynomials:\footnote{Remark: In the paper \cite[Theorem 2.2(a)]{armstrong} we have shown that for any integers $m,n$ the polynomial $(\Box[m,n])_q$, not divided by $[a]_q$, satisfies cyclic sieving for the action of the group $\ZZ/(a-1)\ZZ$ on the set of lattice points $\L \cap\Box[m,n]$ by cyclic permutation of coordinates.}
\begin{equation*}
\frac{1}{[a]_q}(\Box[m,n])_q := \frac{1}{[a]_q}\,\sum_{\x\in(\L\cap \Box[m,n])} q^{\Tsf(\x)} \in\ZZ[q].
\end{equation*}
\end{theorem}

\begin{proof}
For any subset $S\subseteq\RR^{a-1}$ we write $(S)_q$ for the sum of $q^{\Tsf(\x)}$ over the integer points $\L\cap S$, assuming that this sum is finite. First we will compute $(\Box\cap n\Delta)_q$ using an inclusion-exclusion argument similar to the proof of Theorem \ref{thm:johnson2}. For each subset $I\subseteq\{1,\ldots,a-1\}$ we again let $\y_I$ be the indicator vector with $i$-th coordinate equal to $1$ when $i\in I$ and $0$ when $i\not\in I$. Define the cone $K_I=a\y_I+K$ and the simplex $\Delta_I=K_I\cap  n\Delta$, where $K$ is the positive cone. Then we have
\begin{align*}
\Delta_I &= \{a\y_I+\x:  x_i\geq 0 \text{ and } \sum_i (ay_i+x_i)\leq n\}\\
&= \{a\y_I +\x:  x_i\geq 0 \text{ and } \sum_i x_i\leq n-a\#I \}\\
&= a\y_I +(n-a\#I)\Delta,
\end{align*}
and hence
\begin{equation*}
(\Delta_I)_q = q^{a\#I} ((n-a\#I)\Delta)_q= q^{a\#I} \sqbinom{a-1+n-a\#I}{a-1}_q.
\end{equation*}
Since $\Delta_I=\bigcap_{i\in I} \Delta_{\{i\}}$ and since $\Box\cap n\Delta$ equals $\Delta_\emptyset=n\Delta$ minus $\bigcup_{I\neq \emptyset}\Delta_I$ we have
\begin{align*}
(\Box\cap n\Delta)_q &= \sum_{I\subseteq \{1,\ldots,a-1\}} (-1)^{\#I} (\Delta_I)_q \\
&=\sum_{I\subseteq \{1,\ldots,a-1\}} (-1)^{\#I} q^{a\#I} \sqbinom{a-1+n-a\#I}{a-1}_q.
\end{align*}
For any two integers $m\leq n$ this implies 
\begin{align*}
(\Box[m,n])_q &= (\Box\cap n\Delta)_q - (\Box\cap (m-1)\Delta)_q\\
&= \sum_{I\subseteq \{1,\ldots,a-1\}} (-1)^{\#I} q^{a\#I} \left(\sqbinom{a-1+n-a\#I}{a-1}_q-\sqbinom{a-1+m-1-a\#I}{a-1}_q\right).
\end{align*}
If $\gcd(a,n)=\gcd(a,m-1)=1$ then we also have $\gcd(a,n-a\#I)=\gcd(a,m-1-a\#I)=1$ for each $I$. Hence dividing both sides by $[a]_q$ gives an alternating sum of rational $q$-Catalan numbers
\begin{equation*}
\frac{1}{[a]_q} (\Box[m,n])_q =\sum_{I\subseteq \{1,\ldots,a-1\}} (-1)^{\#I} q^{a\#I} \left(\Cat(a,n-a\#I)_q-\Cat(a,m-1-a\#I)_q\right),
\end{equation*}
which is a polynomial in $\ZZ[q]$. 
\end{proof}

The most interesting case occurs when $m$ and $n$ are as close together as possible. We define a special notation for these minimal slices of the fundamental box.

\begin{definition}\label{def:germs}
For any integer $c$ with $\gcd(a,c)=1$, let $c'$ be the largest integer less than $c$ satisfying $\gcd(a,c')=1$. Then we define the {\em $q$-Catalan germ}
\begin{equation*}
\Cat((a;c))_q := \frac{1}{[a]_q} (\Box[c'+1,c])_q.
\end{equation*}
Note that $\Cat((a;c))_q=0$ unless $1\leq c\leq (a-1)^2$.
\end{definition}

The name ``Catalan germ'' is inspired by the second formula in the following proposition.

\begin{proposition}\label{prop:germs}
Let $\gcd(a,b)=1$. Then we have
\begin{align*}
[a]_{q^2} [a]_{q^3} \cdots [a]_{q^{a-1}} &= \sum_{\gcd(a,c)=1} \Cat((a;c))_q,\\
\Cat(a,b)_q &= \sum_{\gcd(a,c)=1} \Cat((a;c))_q \sqbinom{a-1+\lfloor (b-c)/a \rfloor}{a-1}_{q^a},
\end{align*}
where the summands are zero unless $1\leq c\leq (a-1)^2$.
\end{proposition}

\begin{proof}
To prove the first formula we note that the slices $\Box[c'+1,c]$ cover all of the integer points in the fundamental Box. Thus we have
\begin{align*}
\sum_{\gcd(a,c)=1} \Cat((a;c))_q &= (\Box)_q / [a]_q \\
&= \frac{1}{[a]_q} \sum_{ \x\in\{0,\ldots,a-1\}^{a-1}} q^{x_1+2x_2+\cdots+(a-1)x_{a-1}}\\
&= \frac{1}{[a]_q} \sum_{x_1=0}^{a-1} q^{x_1} \sum_{x_2=0}^{a-1} (q^2)^{x_2} \cdots \sum_{x_{a-1}=0}^{a-1} (q^{a-1})^{x_{a-1}} \\
&=  \frac{1}{[a]_q}\prod_{i=1}^{a-1} [a]_{q^i}.
\end{align*}
To prove the second formula we use the fact that $\Box$ is a fundamental domain for $\L/a\L$ to express $b\Delta$ as the disjoint union of the sets $(a\x+\Box)\cap b\Delta$ for $\x\in \NN^{a-1}$, only finitely many of which are non-empty. We note that
\begin{equation*}
(a\x+\Box)\cap b\Delta = a\x + \Box[0,b-a\sum_i x_i],
\end{equation*}
hence
\begin{align*}
\frac{1}{[a]_q} ((a\x+\Box)\cap b\Delta)_q &= q^{a\Tsf(\x)}\frac{1}{[a]_q} (\Box[0,b-a\sum_i x_i])_q\\
&= q^{a\Tsf(\x)} \sum_{\substack{ 1\leq c\leq b-a\sum_i x_i \\ \gcd(a,c)=1}} \Cat((a;c))_q.
\end{align*}
Then summing over the partition of $b\Delta$ gives
\begin{align*}
\Cat(a,b)_q &= (b\Delta)_q/[a]_q\\
&= \sum_{\x\in\NN^{a-1}} q^{a\Tsf(\x)}\frac{1}{[a]_q} ((a\x+\Box)\cap b\Delta)_q\\
&= \sum_{\x\in\NN^{a-1}} q^{a\Tsf(\x)} \sum_{\substack{ 1\leq c\leq b-a\sum_i x_i \\ \gcd(a,c)=1}} \Cat((a;c))_q \\
&= \sum_{\gcd(a,c)=1} \Cat((a;c))_q \sum_{\substack{ \x\in \NN^{a-1} \\ \sum_i x_i \leq \lfloor (b-c)/a \rfloor}} (q^a)^{\Tsf(\x)} \\
&=\sum_{\gcd(a,c)=1 } \Cat((a;c))_q \sqbinom{a-1+\lfloor (b-c)/a \rfloor}{a-1}_{q^a}.
\end{align*}
\end{proof}

When $a=3$ the first formula says that
\begin{equation*}
\Cat((3;1))_q + \Cat((3;2))_q+\Cat((3;4))_q = 1+q^2+q^4 = [3]_{q^2}
\end{equation*}
and the second formula is equivalent to the example after Theorem \ref{thm:johnson1}.

\begin{table}
\begin{equation*}
\begin{array}{c|c|l}
\downstrut
a & c & \Cat((a;c))_q\\
\hline
3 & 1 & 1\\
3 & 2 & q^2 \\
3 & 4 & q^4 \\
\hline
4 & 1 & 1\\
4 & 3 &{q}^{2}+{q}^{3}+{q}^{4}+{q}^{6}\\
4 & 5 &{q}^{5}+{q}^{6}+{q}^{7}+{q}^{8}+{q}^{9}+{q}^{10}\\
4 & 7 &{q}^{9}+{q}^{11}+{q}^{12}+{q}^{13}\\
4 & 9 &{q}^{15}\\
\hline
5 & 1 & 1\\
5 & 2 & {q}^{2}+{q}^{4}\\
5 & 3& {q}^{3}+{q}^{5}+{q}^{6}+{q}^{8}\\
5 & 4 & {q}^{4}+{q}^{6}+{q}^{7}+{q}^{8}+{q}^{9}+{q}^{10}+{q}^{12}\\
5 & 6 &
{q}^{6}+{q}^{7}+2\,{q}^{8}+2\,{q}^{9}+2\,{q}^{10}+3\,{q}^{11}+3\,{q}^{12}+2\,{q}^{13}+3\,{q}^{14}+{q}^{15}+2\,{q}^{16}+{q}^{17}+{q}^{18}
\\
5 & 7 &
{q}^{10}+{q}^{11}+{q}^{12}+2\,{q}^{13}+{q}^{14}+2\,{q}^{15}+2\,{q}^{16}+{q}^{17}+2\,{q}^{18}+{q}^{19}+{q}^{20}+{q}^{21}
\\
5 & 8 &
{q}^{12}+2\,{q}^{14}+{q}^{15}+2\,{q}^{16}+2\,{q}^{17}+{q}^{18}+2\,{q}^{19}+2\,{q}^{20}+{q}^{21}+2\,{q}^{22}+{q}^{24}
\\
5 & 9 &
{q}^{15}+{q}^{16}+{q}^{17}+2\,{q}^{18}+{q}^{19}+2\,{q}^{20}+2\,{q}^{21}+{q}^{22}+2\,{q}^{23}+{q}^{24}+{q}^{25}+{q}^{26}
\\
5 & 11
&{q}^{18}+{q}^{19}+2\,{q}^{20}+{q}^{21}+3\,{q}^{22}+2\,{q}^{23}+3\,{q}^{24}+3\,{q}^{25}+2\,{q}^{26}+2\,{q}^{27}+2\,{q}^{28}+{q}^{29}+{q}^{30}
 \\
5 &  12&
{q}^{24}+{q}^{26}+{q}^{27}+{q}^{28}+{q}^{29}+{q}^{30}+{q}^{32}\\
5 &  13& {q}^{28}+{q}^{30}+{q}^{31}+{q}^{33}\\
5 & 14 & {q}^{32}+{q}^{34}\\
5 & 16 & {q}^{36}
\end{array}
\end{equation*}
\caption{List of $q$-Catalan germs for $3\leq a\leq 5$.}
\label{tab:qcat_germs}
\end{table}

Theorem \ref{thm:germs} implies that $q$-Catalan germs are polynomials with integer coefficients. In Table  \ref{tab:qcat_germs} we list the polynomials $\Cat((a;c))_q$ for all $3\leq a\leq 5$. One can observe from the table that all of the coefficients are non-negative. We state this as a conjecture.

\begin{conjecture}\label{conj:positive_germs}
For any integer $c$ with $1\leq c\leq (a-1)^2$ and $\gcd(a,c)=1$, the corresponding $q$-Catalan germ has non-negative coefficients:
\begin{equation*}
\Cat((a;c))_q \in \NN[q].
\end{equation*}
\end{conjecture}
This conjecture implies that the more general polynomials from Theorem \ref{thm:germs} also have non-negative coefficients because
\begin{equation*}
\frac{1}{[a]_q} (\Box[m,n])_q = \sum_{\substack{m\leq c\leq n \\ \gcd(a,c)=1}} \Cat((a;c))_q.
\end{equation*}

The $q$-Catalan germs are relatively inexpensive to compute, thus we have tested this conjecture up to $a=20$. In the next section we will show that the positivity of $q$-Catalan germs actually implies the existence of Johnson statistics (Conjecture \ref{conj:johnson_statistics}). Then in Section \ref{sec:ribbons} we discuss the small examples $a\in\{3,4,5\}$.

\section{A greedy algorithm}\label{sec:greedy}

For any set $S\subseteq\L\otimes\RR = \RR^{a-1}$ we again use the notation  $(S)_q$ for the tilted height generating function over lattice points in this set:
\begin{equation*}
(S)_q = \sum_{\x\in (\L\cap S)} q^{\Tsf(\x)}.
\end{equation*}
We begin with the simple observation that 
\begin{equation*}
(\Box)_q = \sum_{\x\in(\L\cap\Box)} q^{x_1+2x_2+\cdots+(a-1)x_{a-1}} =[a]_q [a]_{q^2} \cdots [a]_{q^{a-1}}.
\end{equation*} 
If $\Jsf:\R\to \ZZ$ is some hypothetical Johnson statistic then we write $(S)_\Jsf$ for the sum of $q^{\Jsf(\x)}$ over points $\x$ in the intersection of $S$ with the root lattice:
\begin{equation*}
(S)_\Jsf = \sum_{\x\in (\R\cap S)} q^{\Jsf(\x)}.
\end{equation*}
From Theorem \ref{thm:johnson2} we have $(\Box)_\Jsf = (\Box)_q / [a]_q$. This suggests that we should focus on the full weight lattice and not just the root lattice when looking for Johnson statistics.

\begin{definition}\label{def:standard}
Consider a set of lattice points $S\subseteq\L$. If $\#S=a$ and if the points of $S$ have consecutive tilted heights, i.e., if $(S)_q = q^k + q^{k+1} + \cdots + q^{k+a-1} = q^k [a]_q$ for some integer $k$, then we say that $S$ is a {\em standard set}.
\end{definition}

For example, the vertices of the fundamental alcove $\{\omega_0,\omega_1,\ldots,\omega_{a-1}\}$ are a standard set because $\Tsf(\omega_i)=i$. More generally, the points of any alcove are a standard set. However, not every standard set is an alcove. For example, for any point $\y\in\L$ the ``line segment'' $\{\y,\y+\omega_1,\y+2\omega_1,\ldots,\y+(a-1)\omega_1\}$ is a standard set because $\Tsf(\y+i\omega_1)=\Tsf(\y)+i$.

Note that each standard set contains a unique point of the root lattice, i.e., a unique point with height $\equiv 0$ mod $a$. Suppose that we can find a partition of the $a^{a-1}$ points of $\L\cap\Box$ into $a^{a-2}$ standard sets:
\begin{equation*}
\L\cap\Box = \coprod_{i=1}^{a^{a-2}} S_i.
\end{equation*}
Let's say $(S_i)_q = q^{k_i} [a]_q$ for some integers $k_i$. If $\x_i\in S_i$ is the unique point of $\R\cap S_i$ then we define a statistic $\Jsf:(\R\cap\Box)\to\ZZ$ by $\Jsf(\x_i):=k_i$. That is, we  define $\Jsf(\x_i)=\Tsf(\y_i)$, where $\y_i\in S_i$ is the lowest point in this ribbon. Note that this statistic satisfies
\begin{equation*}
(\Box)_q = \sum_{i=1}^{a^{a-2}} (S_i)_q = \sum_{i=1}^{a^{a-2}} q^{k_i} [a]_q= [a]_q \sum_{i=1}^{a^{a-2}} q^{k_i} = [a]_q \sum_{\x\in (\R\cap\Box)} q^{\Jsf(\x)} =[a]_q (\Box)_\Jsf.
\end{equation*}
It follows that $(\Box)_\Jsf= [a]_{q^2} \cdots [a]_{q^{a-1}}$, which is one of the desired properties of Johnson statistics. In fact, it is quite easy to find partitions of the fundamental box into standard sets. For example, for each choice of the coordinates $\x'=(x_2,\ldots,x_{a-1})$ satisfying $0\leq x_i\leq a-1$ there is a standard set
\begin{equation*}
S_{\x'} = \{(x_1,\x'): 0\leq x_1\leq a-1\}
\end{equation*}
satisfying
\begin{equation*}
(S_{\x'})_q = \sum_{x_1=0}^{a-1} q^{x_1+2x_2+\cdots+(a-1)x_{a-1}} = q^{2x_2+\cdots+(a-1)x_{a-1}} [a]_q.
\end{equation*}
Thus the function $\Jsf:(\R\cap\Box)\to\ZZ$ defined by
\begin{equation*}
\Jsf(x_1,x_2,\ldots,x_{a-1}) := 2x_2+\cdots+ (a-1)x_{a-1}
\end{equation*}
satisfies $(\Box)_\Jsf =[a]_{q^2} \cdots [a]_{q^{a-1}}$. And we can easily extend this statistic to all of $\R$ by periodicity, setting $\Jsf(\x+a\y):=\Jsf(\x)+a\Tsf(\y)$ for all $\y\in\L$.

But it is much more difficult to find a statistic $\Jsf:\R\to\ZZ$ satisfying the Catalan property of Johnson statistics:
\begin{equation*}
(b\Delta)_\Jsf = \Cat(a,b)_q \text{ \,\,for all integers $b\geq 1$ such that $\gcd(a,b)=1$.}
\end{equation*}
To find such a statistic we should look for a standard partition of $\L\cap\Box$ in which the standard sets do not cross hyperplanes of the form $x_1+\cdots+x_{a-1}=c+1/2$ for $\gcd(a,c)=1$. In other words, we should look for a standard partition of each set $\L\cap\Box[c'+1,c]$ where $\gcd(a,c)=1$ and where $c'<c$ is the largest integer less than $c$ satisfying $\gcd(a,c')=1$. Recall we have shown that $\Cat((a;c))_q=(\Box[c'+1,c])_q/[a]_q$ is in $\ZZ[q]$. 

\begin{conjecture}\label{conj:standard}
For any integer $c$ with $1\leq c\leq (a-1)^2$ and $\gcd(a,c)=1$, the set $\L\cap\Box[c'+1,c]$ has a standard partition.
\end{conjecture}

I claim that this conjecture implies Conjecture \ref{conj:johnson_statistics}.

\begin{proof}
If $\L\cap\Box[c'+1,c]$ has a standard partition then, as in the previous proof, we obtain a statistic $\Jsf:(
\R\cap\Box[c'+1,c])\to\ZZ$ with the property
\begin{equation*}
(\Box[c'+1,c])_\Jsf = (\Box[c'+1,c])_q / [a]_q =\Cat((a;c))_q.
\end{equation*}
But $\R\cap\Box$ is the disjoint union of the sets $\R\cap\Box[c'+1,c]$ over $1\leq c\leq (a-1)^2$ with $\gcd(a,c)=1$. Let us concatenate these to define a statistic $\Jsf:(\R\cap\Box)\to\ZZ$ and then extend $\Jsf$ to all of $\R$ by setting $\Jsf(\x+a\y)=\Jsf(\x)+a\Tsf(\y)$ for all $\y\in\L$. I claim that this $\Jsf$ is a Johnson statistic, and to show this we only need to check the Catalan property. That is, for any $\gcd(a,b)=1$ we need to show that $(b\Delta)_\Jsf = \Cat(a,b)_q$. 

As an intermediate step, we note that for any $\gcd(a,b)=1$ we have
\begin{align*}
(\Box[0,b])_q /[a]_q &=  \frac{1}{[a]_q}\sum_{\substack{1\leq c\leq b\\ \gcd(a,c)=1}} (\Box[c'+1,c])_q\\
&=  \sum_{\substack{1\leq c\leq b\\ \gcd(a,c)=1}} (\Box[c'+1,c])_q/[a]_q\\
&=  \sum_{\substack{1\leq c\leq b\\ \gcd(a,c)=1}}(\Box[c'+1,c])_\Jsf \\
&=  \biggl(\,\coprod_{\substack{1\leq c\leq b\\ \gcd(a,c)=1}}\Box[c'+1,c]\biggr)_\Jsf \\
&= (\Box[0,b])_\Jsf.
\end{align*}
Then we recall the decomposition of $b\Delta$ from the proof of \ref{prop:germs}:
\begin{align*}
\Cat(a,b)_q &= \frac{1}{[a]_q} (b\Delta)_q\\
&= \frac{1}{[a]_q} \sum_{\y\in\NN^{a-1}} (a\y+\Box[0,b-a\sum_i y_i])_q\\
&= \frac{1}{[a]_q} \sum_{\y\in\NN^{a-1}} q^{a\Tsf(\y)} (\Box[0,b-a\sum_i y_i])_q\\
&= \sum_{\y\in\NN^{a-1}} q^{a\Tsf(\y)} \frac{1}{[a]_q}(\Box[0,b-a\sum_i y_i])_q\\
&= \sum_{\y\in\NN^{a-1}} q^{a\Tsf(\y)} \sum_{\x\in (\R\cap \Box[0,b-a\sum_i y_i])} q^{\Jsf(\x)}\\
&= \sum_{\y\in\NN^{a-1}} \, \sum_{\x\in (\R\cap \Box[0,b-a\sum_i y_i])} q^{\Jsf(\x)+a\Tsf(\y)}\\
&= \sum_{\y\in\NN^{a-1}} \, \sum_{\x\in (\R\cap \Box[0,b-a\sum_i y_i])} q^{\Jsf(\x+a\y)}\\
&= \sum_{\x \in (\R\cap b\Delta)} q^{\Jsf(\x)}\\
&= (b\Delta)_\Jsf.
\end{align*}
\end{proof}

Finally, we show that Conjectures \ref{conj:positive_germs} and \ref{conj:standard} are equivalent, from which it follows that Conjecture \ref{conj:standard}, and hence also Conjecture \ref{conj:johnson_statistics}, is true for all $a\leq 20$.

\begin{proof}
Fix positive integers $a,c$ with $1\leq c\leq (a-1)^2$ and $\gcd(a,c)=1$. If $\L\cap\Box[c'+1,c]$ has a standard partition then as in the previous proof then we have
\begin{equation*}
\Cat((a;c))_q = \sum_{\x\in(\R\cap\Box[c'+1,c])} q^{\Jsf(\x)},
\end{equation*}
which implies that $\Cat((a;c))_q$ has non-negative coefficients. Conversely, suppose that $\Cat((a;c))_q$ has non-negative coefficients. Let $Y=\L\cap\Box[c'+1,c]$ and consider the polynomial 
\begin{equation*} 
F(q):=(Y)_q = \sum_{\y\in  Y} q^{\Tsf(\y)}.
\end{equation*}
By assumption we have $F(q)/[a]_q = \Cat((a;c))_q\in\NN[q]$ so that $F(q)=q^k [a]_a G(q)$ for some integer $k\geq 0$ and polynomial $G(q)\in \NN[q]$ with non-zero constant term. If $g_0\geq 1$ is the constant term of $G(q)$ then $F(q) = g_0q^k+g_0q^{k+1}+\cdots+g_0 q^{k+a-1} + q^{k'}G'(q)$ for some integer $k'\geq k+1$ and polynomial $G'(q)\in\NN[q]$ with non-zero constant term. Since the coefficients of $q^k,q^{k+1},\ldots,q^{k+a-1}$ in $F(q)$ are all $\geq g_0$, it follows that we can find $g_0$ standard sets of points in $Y$ with heights $k,k+1,\ldots,k+a-1$.\footnote{In the first step we will have $g_0=1$, but we allow the general case for the purpose of induction.} Let $Y'$ denote the subset of $Y$ with all of these standard sets deleted. The new generating function is
\begin{equation*}
F'(q):= (Y')_q = \sum_{\y\in  Y'} q^{\Tsf(\y)} = q^{k'}G'(q).
\end{equation*}
But $F'(q)/[a]_q$ still has non-negative coefficients, so we may proceed by induction until we have partitioned $Y$ into standard sets.
\end{proof}

We can rephrase this proof as an algorithm. Consider any finite set $S\subseteq\L$ with the property that $(S)_q/[a]_q$ is a polynomial with non-negative integer coefficients. Then we can partition $S$ into standard sets as follows:
\begin{itemize}
\item Let $\x_0\in S$ be any point of minimal tilted height.
\item Greedily choose points $\x_1,\ldots,\x_{a-1}\in S$ satisfying $\Tsf(\x_i)=\Tsf(\x_0)+i$.
\item Replace $S$ by $S\setminus\{\x_0,\x_1,\ldots,\x_{a-1}\}$.
\item Repeat.
\end{itemize}
The above proof shows that this procedure never gets stuck.

\section{Ribbon partitions}\label{sec:ribbons}

The standard partitions produced by the greedy algorithm can be quite random. In order to prove existence in general it is reasonable to consider the partial order on the weight lattice.

\begin{definition}\label{def:ribbons}
Recall from Section \ref{sec:tilted} that for integer points $\x,\y\in \L$ we write $\x\leq \y$ if and only if $x_i+x_{i+1}+\cdots+x_{a-1}\leq y_i+y_{i+1}+\cdots+y_{a-1}$ for all $1\leq i\leq a-1$. A saturated chain of length $a$ in this partial order is called a {\em ribbon}. Note that each ribbon is a standard set because the tilted height is a rank function for the tilted partial order.
\end{definition}

The following strengthening of Conjecture \ref{conj:standard} is based on rather limited evidence from the small cases $a\in\{3,4,5\}$, but we still state it as a conjecture.

\begin{conjecture}\label{conj:ribbon}
Choose any integer $1\leq c\leq (a-1)^2$ with $\gcd(a,c)=1$, and let $c'$ be the largest integer less than $c$ satisfying $\gcd(a,c')=1$. Then the tilted partial order on the set of lattice points $\L\cap\Box[c'+1,c]$ has a partition into disjoint ribbons.
\end{conjecture}

It may be helpful to translate this conjecture into the language of Young diagrams. Given a point $\x=x_1\omega_1+\cdots+x_{a-1}\omega_{a-1}$ expressed in weight coordinates, we can write $\x=y_1\omega_1+y_2(\omega_2-\omega_1)+\cdots+y_{a-1}(\omega_{a-1}-\omega_{a-2})$ in tilted coordinates. Recall that the tilted partial order is just the componentwise order on tilted coordinates, which we can view as inclusion of Young diagrams. For integers $m\leq n$ the integer points in the slice $\Box[m,n]$ expressed in tilted coordinates are
\begin{multline*}
\L\cap \Box[m,n] = \{(y_1,\ldots,y_{a-1})\in\ZZ^{a-1}: y_1\geq \cdots\geq y_{a-1}\geq 0\\
\text{ with }m\leq y_1\leq n  \text{ and } 0\leq y_i-y_{i+1}\leq a-1
\text{ for } 1\leq i\leq a-2\}.
\end{multline*}
These are the Young diagrams with at most $a$ columns, largest row length between $m$ and $n$, and difference between consecutive row lengths $\leq a-1$. If $\gcd(a,m-1)=\gcd(a,n)=1$ then Conjecture \ref{conj:ribbon} implies that the inclusion order on these Young diagrams has a ribbon partition, i.e., a partition into disjoint saturated chains of length $a$.

\begin{figure}
\begin{center}
\includegraphics[scale=0.8]{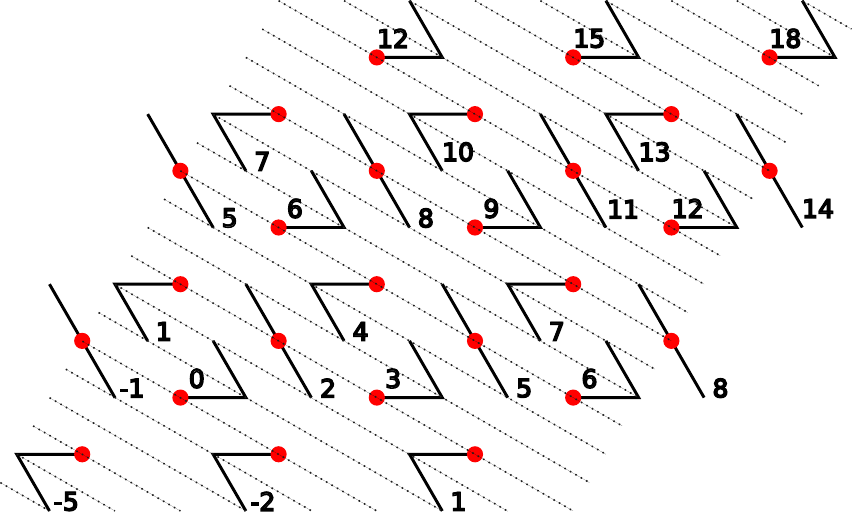}
\end{center}
\caption{The unique ribbon partition when $a=3$.}
\label{fig:ribbon_partition_a2}
\end{figure}

In the remainder of this section we show that Conjecture \ref{conj:ribbon} holds for $a\in\{3,4,5\}$. In the case $a=3$ there are only three integers $1\leq c\leq (3-1)^2$ with $\gcd(3,b)=1$; namely, $c\in\{1,2,4\}$. The corresponding sets of points are
\begin{align*}
\L\cap \Box[0,1] &= \{\langle 0,0\rangle ,(1,0),(0,1)\},\\
\L\cap \Box[2,2]  &= \{(2,0),\langle 1,1\rangle ,(0,2)\},\\
\L\cap \Box[3,4]  &= \{(2,1),(1,2),\langle 2,2\rangle \}.
\end{align*}
Each of these posets is itself a ribbon, so there exists a unique ribbon partition. Note that each ribbon contains a unique point of the root lattice, which we write using angle brackets. To obtain a Johnson statistic $\Jsf:\R\to\ZZ$ from a ribbon partition we define $\Jsf(\x)=\Tsf(\y)$ where $\y\in\L$ is the lowest point in the ribbon containing $\x$. That is, we define
\begin{align*}
\Jsf(0,0) &:= \Tsf(0,0) = 0,\\
\Jsf(1,1) &:= \Tsf(2,0) = 2,\\
\Jsf(2,2) &:= \Tsf(2,1) = 4.
\end{align*}
And then we extend $\Jsf$ to the full root lattice by periodicity, setting $\Jsf(\x+a\y):=\Jsf(\x)+a\Tsf(\y)$ for all $\x\in(\R\cap\Box)$ and $\y\in\L$. This is precisely the example presented in Section \ref{sec:johnson}. Figure \ref{fig:ribbon_partition_a2} shows this unique ribbon partition, with points of the root lattice shown in red. Compare to Figure \ref{fig:johnson_a2}.

\begin{table}[h]
\begin{equation*}
\begin{array}{cccc|c}
&\hfill\text{ribbon}&\text{points}\hfill&& \Jsf \\
\hline
\upstrut\downstrut
(0,0,1) & (0,1,0) & (1,0,0) & \langle 0,0,0\rangle  & 0\\
\hline
\upstrut
(2,0,1) & \langle 2,1,0\rangle & (3,0,0) & (2,0,0) & 2\\
(1,1,1) &(0,1,1) & \langle 1,0,1\rangle & (1,1,0) & 3\\
(0,2,1) &(0,3,0) &(1,2,0)&\langle 0,2,0\rangle & 4\\
\downstrut
(0,0,3) &\langle 0,1,2\rangle &(1,0,2)& (0,0,2) & 6\\
\hline
\upstrut
\langle 3,1,1\rangle &(2,1,1)&(3,0,1)& (3,1,0) & 5\\
(2,2,1) &\langle 2,3,0\rangle &(3,2,0)& (2,2,0) & 6\\
(1,3,1) &(0,3,1)&\langle 1,2,1\rangle & (1,3,0) & 7\\
(2,0,3) &(2,1,2)&(3,0,2) & \langle 2,0,2\rangle & 8\\
\langle 1,1,3\rangle  &(0,1,3)&(1,0,3)& (1,1,2) & 9\\
\downstrut
(0,2,3) &\langle 0,3,2\rangle &(1,2,2)& (0,2,2) & 10\\
\hline
\upstrut
\langle 3,3,1\rangle &(2,3,1)&(3,2,1)& (3,3,0) & 9\\
(3,1,3) &(2,1,3)&\langle 3,0,3\rangle & (3,1,2) & 11\\
(2,2,3) &(2,3,2)&(3,2,2)& \langle 2,2,2\rangle & 12\\
\downstrut
\langle 1,3,3\rangle &(0,3,3) &(1,2,3)& (1,3,2) & 13\\
\hline
\upstrut
(3,3,3) &(2,3,3)&\langle 3,2,3\rangle& (3,3,2) & 15
\end{array}
\end{equation*}
\caption{The remarkable ribbon partition when $a=4$.}
\label{tab:ribbon_partition_a3}
\end{table}

\begin{figure}[h]
\begin{center}
\includegraphics[scale=0.5]{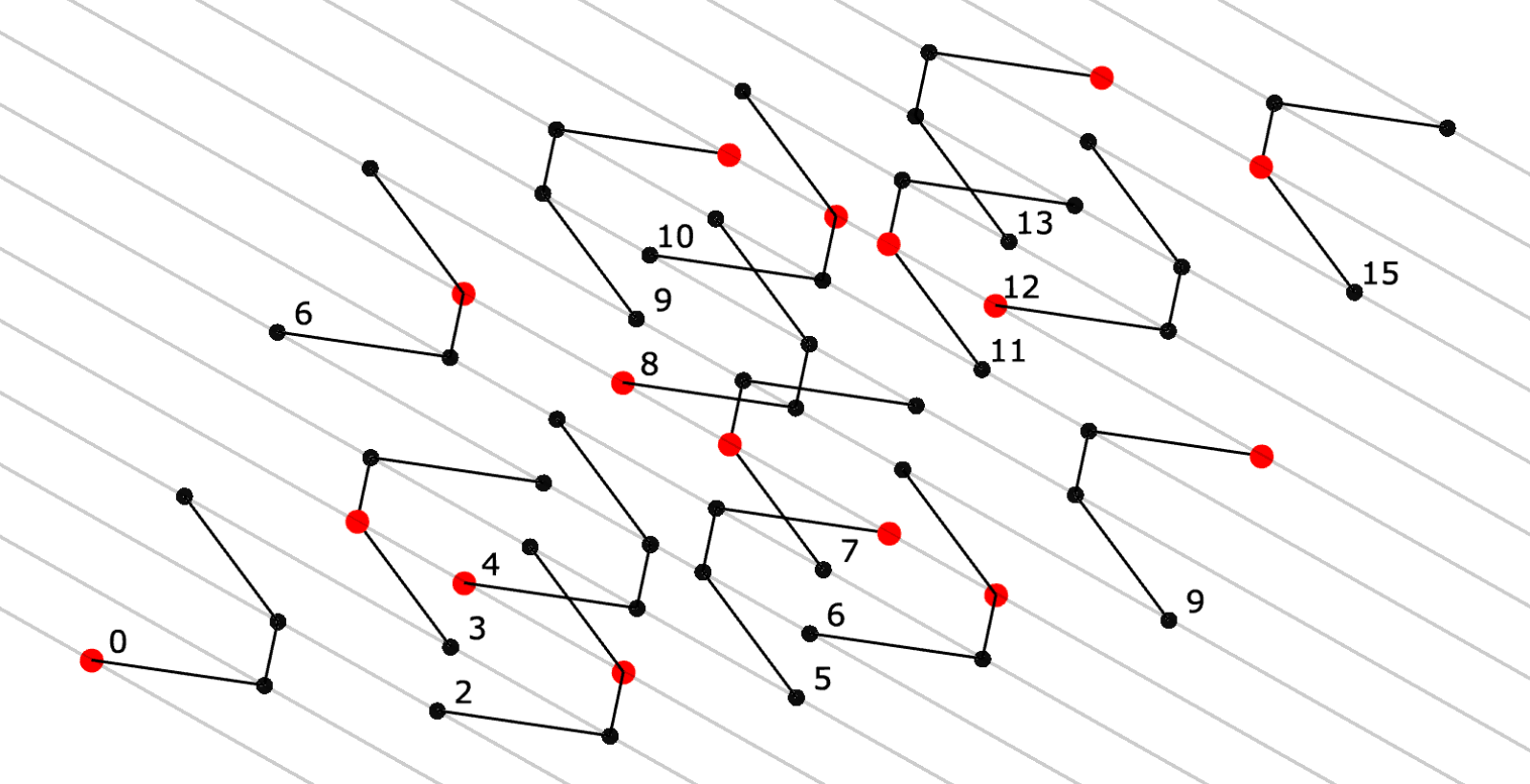}
\end{center}
\caption{The remarkable ribbon partition when $a=4$.}
\label{fig:ribbon_partition_a3}
\end{figure}

Next we consider $a=4$. In this case there exists a ribbon partition with remarkably nice geometric properties. Namely, we consider the infinite hyperplane arrangement consisting of hyperplanes of the form $x_i=2k-1/2$ for $i\in\{1,2,3\}$ and $k\in\ZZ$, and the hyperplanes $x_1+x_2+x_3=c+1/2$ for $\gcd(4,c)=1$, i.e., for odd integers $c\in\ZZ$. It turns out that each connected component of this hyperplane arrangement contains exactly $4$ lattice points, which are the vertices of an alcove, and hence form a ribbon. By periodicity, it is enough to list these ribbons for the $4^{4-1}=64$ points in the fundamental box, which we display in Table \ref{tab:ribbon_partition_a3}. Each row of the table contains the points of a ribbon, with the root lattice point written in angle brackets. For each point $\x\in(\R\cap\Box)$ we define $\Jsf(\x):=\Tsf(\y)$, where $\y$ is the lowest point in the ribbon (the rightmost point in the row). The value of $\Jsf$ is shown in the rightmost column. For example, we have
\begin{equation*}
\Jsf(1,1,3):=\Tsf(1,1,2)=1\cdot 1+2\cdot 1+3\cdot 2=9.
\end{equation*}
Horizontal lines in the table separate the points in the different slices of $\Box$, corresponding to the $q$-Catalan germs. For example, the slice $\Box[6,7]$ contains $16$ points of $\L$ and $4$ points of $\R$, corresponding to the germ
\begin{equation*}
\Cat((4;7))_q = \frac{1}{[4]_q}\,\sum_{\y\in (\L\cap\Box[6,7])} q^{\Tsf(\y)} = \sum_{\x \in (\R\cap\Box[6,7])} q^{\Jsf(\x)} = q^{9}+q^{11}+q^{12}+q^{13}.
\end{equation*}
Compare to Table \ref{tab:qcat_germs}. In Figure \ref{fig:ribbon_partition_a3} we show a projection of this ribbon partition for points in the fundamental box. The diagonal lines connect points with the same tilted height. The points of the root lattice are drawn in red, and one can see that they have tilted heights congruent to $0$ mod $4$, as expected.

\begin{figure}
\begin{center}
\includegraphics[scale=0.8]{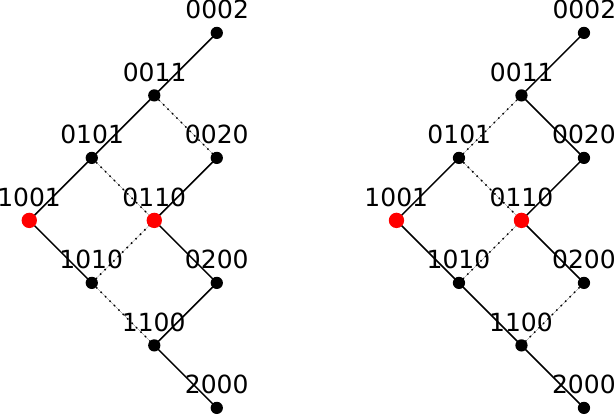}
\end{center}
\vspace{.1in}
\caption{The poset $\L\cap\Box[2,2]$ for $a=5$ has two ribbon partitions.}
\label{fig:box_22}
\end{figure}

After the success of $a=3$ and $a=4$, the case $a=5$ is a bit disappointing. In this case the poset $\L\cap\Box[2,2]$, corresponding to the germ $\Cat((5;2))_q=q^2+q^4$, has exactly two ribbon partitions, and we see no good reason to prefer one over the other. Figure \ref{fig:box_22} shows these two partitions. The red points $(1,0,0,1)$ and $(0,1,1,0)$ are in the root lattice. The lowest points in the corresponding ribbons are either $(1,0,1,0)$ and $(2,0,0,0)$, or $(2,0,0,0)$ and $(0,2,0,0)$, respectively. The first choice yields Johnson statistic $\Jsf(1,0,0,1):=\Tsf(1,0,1,0)=4$ and $\Jsf(0,1,1,0):=\Tsf(2,0,0,0)=2$, while the second choice yields $\Jsf(1,0,0,1):=\Tsf(2,0,0,0)=2$ and $\Jsf(0,1,1,0):=\Tsf(1,0,1,0)=4$. This suggests that there is no ``canonical'' Johnson statistic when $a=5$.  Nevertheless, we have verified by computation that at least one ribbon partition does exist.\footnote{We thank an anonymous referee for independent confirmation of this fact.} Unfortunately a greedy algorithm does not work here so the computations are more expensive. We found the case $a=6$ out of reach.

For the reader's convenience we summarize the dependencies among the conjectures thus far:
\begin{align*}
\left\{\begin{array}{c}\text{\ref{conj:ribbon} ribbon}\\\text{partitions}\end{array}\right\} &\Longrightarrow 
\left\{\begin{array}{c}\text{\ref{conj:standard} standard}\\\text{partitions}\end{array}\right\} \\&\Longleftrightarrow \left\{\begin{array}{c}\text{\ref{conj:positive_germs} positive}\\\text{germs}\end{array}\right\} \\
&\Longrightarrow \left\{\begin{array}{c}\text{\ref{conj:johnson_statistics} Johnson}\\\text{statistics}\end{array}\right\}\\
& \Longrightarrow 
\left\{\begin{array}{c}\text{\ref{conj:monotone} monotone}\\\text{$q$-Catalan}\end{array}\right\}.
\end{align*}
We end this section with a few words about a possible method of proof. Sylvester's proof of unimodality for $q$-binomial coefficients \cite{sylvester} implicitly used an action of the group $SL_2$ on spaces of invariants. Stanley \cite{stanley_paper} generalized this idea to other posets arising from Lie groups, and then Proctor \cite{proctor} boiled down the argument to its bare essentials.

We briefly describe Proctor's approach, translated into our language. Let $L(a-1,b)$ be the poset of Young diagrams fitting in an $(a-1)\times b$ rectangle, which is the same as our tilted partial order on the lattice points $\L\cap b\Delta$. The rank-generating polynomial of this poset is the $q$-binomial coefficient:
\begin{equation*}
(b\Delta)_q = \sum_{\x\in (\L\cap b\Delta)} q^{\Tsf(\x)} = \sqbinom{a-1+b}{a-1}_q.
\end{equation*}
Let $V=V_0\oplus V_1\oplus \cdots\oplus V_{(a-1)b}$ where $V_k$ is the vector space abstractly generated by the points $\x\in (\L\cap b\Delta)$ satisfying $\Tsf(\x)=k$, so that
\begin{equation*}
\sqbinom{a-1+b}{a-1}_q = \sum_k \dim(V_k) q^k.
\end{equation*}
For each point $\x\in V_k$, let $X(\x)\in V_{k+1}$ be the formal sum of points that cover $\x$ in the poset, setting $X(\x)=\0$ if there are no such points. Proctor showed that there exists a ``lowering operator'' $Y:
V\to V$ satisfying $Y(V_k)\subseteq V_{k-1}$, where $H:=XY-YX$ is a scalar multiple of the identity satisfying $HX-XH=2X$ and $HY-YH=-2Y$. It then follows from the representation theory of $SL_2$ that the vector space $V$ decomposes into subspaces $U_i$ with the property that $\dim(V_k\cap U_i)\in\{0,1\}$ and such that each $U_i$ is supported on a consecutive string of ranks centered at rank $(a-1)b/2$. This is enough to show that the sequence of coefficients is unimodal. But this argument does not prove the existence of a symmetric chain decomposition, even though the operators are defined in terms of the partial order.

We can imagine a similar method to produce standard partitions and hence Johnson statistics. Given one of the posets $\L\cap \Box[c'+1,c]$ we let $V=V_0\oplus V_1\oplus\cdots\oplus V_{a-1}$ be the cyclically-graded vector space where $V_k$ is the formal sum of points $\x\in(\L\cap \Box[c'+1,c])$ satisfying $\Tsf(\x)\equiv k$ mod $a$. One might use the poset structure to construct an operator $X:V\to V$ that is nilpotent of order $a$ and satisfies $X(V_k)\subseteq V_{(k+1 \text{ mod } a)}$. Hopefully this operator would have special properties that can be used to decompose the vector space $V$ into ``$SL_2$-strings'' of length $a$. As with the problem of symmetric chain decompositions, however, a dimension argument of this type would not be enough to produce ribbon partitions.

\section{Further questions}

We end this paper with some speculation and suggestions for
future work.

\subsection{Rank-level duality}\label{sec:rank_level} Even though we
have emphasized the asymmetric expression
\begin{equation*}
\Cat(a,b)_q = \frac{1}{[a]_q}\sqbinom{a-1+b}{a-1}_q,
\end{equation*}
it is nevertheless true that $\Cat(a,b)_q=\Cat(b,a)_q$. Let
$\R^{(a)}\subseteq\ZZ^{a-1}$ denote the root lattice of dimension
$a-1$. If we can find Johnson statistics $\Jsf:\R^{(a)}\to\ZZ$
and $\Jsf:\R^{(b)}\to\ZZ$ then for any integers $a,b\geq 1$
with $\gcd(a,b)=1$ there should exist a bijection $\R^{(a)}\cap
b\Delta\leftrightarrow\R^{(b)}\cap a\Delta$ preserving
these statistics. Consider the example $\{a,b\}=\{3,4\}$ with
$\Cat(3,4)_q=\Cat(4,3)_q=1+q^2+q^3+q^4+q^6$ (see Figures \ref{fig:ribbon_partition_a2} and \ref{fig:ribbon_partition_a3}). In this case we have
\begin{align*}
\R^{(3)}\cap 4\Delta &= \{(0,0), (1,1), (3,0), (2,2), (0,3) \},\\
\R^{(4)}\cap 3\Delta &= \{(0,0,0), (2,1,0), (1,0,1), (0,2,0), (0,1,2) \}.
\end{align*}
According to the Johnson statistics for $a=3$ and $a=4$ defined in the previous section, the desired bijection is
\begin{equation*}
\begin{array}{c|c|c}
\downstrut
\R^{(3)}\cap 4\Delta & \R^{(4)}\cap 3\Delta & \Jsf \\
\hline
\upstrut
(0,0) & (0,0,0) & 0\\
(1,1) & (2,1,0) & 2\\
(3,0) & (1,0,1) & 3\\
(2,2) & (0,2,0) & 4\\
(0,3) & (0,1,2) & 6
\end{array}
\end{equation*}
This problem is known to be difficult even in the unweighted case
$q=1$.

\subsection{Non-coprime parameters}\label{sec:non_coprime}

If $\gcd(a,b)\neq 1$ then the expression $\Cat(a,b)=\binom{a-1+b}{a-1}/a$ is not an integer. There are different ways to resolve this issue, depending on which interpretation of $\Cat(a,b)$ we are using. We will look at this from the point of view of lattice points. Recall that the cosets of the root lattice in the weight lattice have the form
  \begin{equation*}
 \R_k=\{(x_1,\ldots,x_{a-1})\in\ZZ^{a-1}:
 x_1+2x_2+\cdots+(a-1)x_{a-1}\equiv k \mod a\}.
 \end{equation*}
If $\gcd(a,b)=1$ then we showed in Section \ref{sec:rational_cat} that $\#(\R_k\cap b\Delta)=\Cat(a,b)$ for any $k$. This is no longer the case when $\gcd(a,b)\neq 1$. The general counting formula was found by Fredman \cite[Formula 4]{fredman}. It uses the {\em Ramanujan sum}
\begin{equation*}
c_d(\ell) = \sum_{s |{\gcd(d,\ell)}} \mu\left(\frac{d}{s}\right) s,
\end{equation*}
where $\mu$ is the number-theoretic M\"obius function. We note that $c_d(\ell)=\mu(d)$ when $\gcd(d,\ell)=1$ and $c_d(\ell)=\phi(d)$ is Euler's phi function when $\gcd(d,\ell)=d$. Furthermore, we have $c_d(\ell+dt)=c_d(\ell)$ for any $t\in\ZZ$. If we define $\Cat(a,b;k):=\#(\R_k \cap b\Delta)$ then the theorem says\footnote{The key point of \cite{fredman} is the observation that $\Cat(a,b;k)=\Cat(b,a;k)$ for any $k$. Elashvili, Jibladze and Pataraia \cite{elashvili} later rediscovered this result and they called it ``Hermite reciprocity". Compare to Section \ref{sec:rank_level}. }
\begin{equation*}
\Cat(a,b;k) = \frac{1}{a+b} \sum_{d|{\gcd(a,b)}}  c_d(k) \binom{a/d+b/d}{a/d}.
\end{equation*}
When $\gcd(a,b)=1$ we note that $\Cat(a,b;k)=\Cat(a,b)$ for any $k$. Reiner, Stanton and White \cite{rsw} later defined the {\em $q$-Ramanujan sum} by
\begin{equation*}
c_d(\ell)_q = \sum_{s |{\gcd(d,\ell)}} \mu\left(\frac{d}{s}\right) [s]_q,
\end{equation*}
and they (implicitly) considered the following expression:\footnote{Their setup has an extra parameter and uses an $m$-tuple $(k_1,\ldots,k_m)$ instead of a pair $(a,b)$.}
\begin{equation*}
\Cat(a,b;k)_q := \frac{1}{[a+b]_q} \sum_{d|{\gcd(a,b)}} c_d(k)_q \binom{a/d+b/d}{a/d}_{q^d}.
\end{equation*}
In \cite[Section 9]{rsw} they worked with a finite field $\FF_q$ and the action of the group $\FF_{q^{a+b}}^\times / \FF_q^\times$ on the set of $a$-dimensional (or $b$-dimensional) subspaces of $\FF_{q^{a+b}}$. They claimed that it is ``a straightforward exercise in M\"obius inversion'' to show that $\Cat(a,b;k)_q$ counts the orbits whose stabilizer order divides $[k]_q$. In \cite[Section 10]{rsw} they proved (Corollary 10.2) that $\Cat(a,b;k)_q$ is in $\ZZ[q]$ and they conjectured (Conjecture 10.3) that the coefficients are non-negative. Based on computer experiments we make the following conjecture generalizing their Conjecture 10.3 and our Conjecture \ref{conj:monotone} (the case $g=1$). We have tested this conjecture for $a\leq 20$ and $b\leq c\leq 80$.

\begin{conjecture}\label{conj:monotone2}
Fix integers $g,k\in\NN$. Then for any integers $a,b,c\geq 1$ satisfying $\gcd(a,b)=\gcd(a,c)=g$ we have
\begin{equation*}
\Cat(a,c;k)_q - \Cat(a,b;k)_q \in\NN[q].
\end{equation*}
\end{conjecture}
When $a=3$ and $g=3$ this conjecture says that we should have
\begin{equation*}
\Cat(3,3t;k)_q - \Cat(3,3s;k)_q \in\NN[q]
\end{equation*}
for all integers $k,s,t$ with $0\leq s\leq t$. We observe that this holds for the examples
\begin{align*}
\Cat(3,3;0)_q &=1 + q + q^{2} + q^{4}, \\
\Cat(3,6;0)_q &=1 + q + q^{2} + 2q^{4} + q^{5} + q^{6} + q^{7} + q^{8} + q^{10}, \\
\Cat(3,6;0)_q-\Cat(3,3;0)_q &=q^{4} + q^{5} + q^{6} + q^{7} + q^{8} + q^{10},
\end{align*}
and
\begin{align*}
\Cat(3,3;1)_q &= q + q^{2} + q^{4},\\
\Cat(3,6;1)_q &= q + q^{2} + 2q^{4} + q^{5} + q^{6} + q^{7} + q^{8} + q^{10},\\
\Cat(3,6;1)_q-\Cat(3,3;1)_q &=q^{4} + q^{5} + q^{6} + q^{7} + q^{8} + q^{10}.
\end{align*}
The Johnson statistics in this paper do not explain these phenomena. For example, the unique Johnson statistic with $a=3$ gives $(3\Delta)_\Jsf=1+q^2+q^3+q^6$, which does not agree with the formula $\Cat(3,3;0)_q=1 + q + q^{2} + q^{4}$. Perhaps there exists a modified Johnson statistic $\Jsf_{g,k}:\R_k\to\ZZ$ for each pair $g,k$ but we have not looked into this.

Dyck paths lead to a {\bf different} generalization of rational Catalan numbers to non-coprime parameters, corresponding to paths in a rectangle that stay weakly above the diagonal. These were originally counted by Bizley \cite{bizley}. Bergeron, Garsia, Leven and Xin \cite{bglx} described $q$ and $t$ statistics for these paths, but it is not easy to see what happens under the substitution $t=1/q$. It would be interesting to find a connection between the two pictures.

\subsection{Other Weyl groups}\label{sec:weyl}

The concept of Catalan numbers and their $q$- and $q,t$-analogues
makes sense for any crystallographic root system. Let
$\R\subseteq\L$ be a root lattice and weight lattice of rank $n$
(we refer to Humphreys~\cite{humphreys} for the basic theory). Let
$\omega_1,\ldots,\omega_n\in\L$ denote the basis of fundamental weights,
and let $C$ be the $n\times n$ Cartan matrix. The basis of simple roots
$\alpha_1,\ldots,\alpha_n\in\R$ (expressed in weight coordinates) are
the columns of the Cartan matrix. Therefore the index of connection
$f=\#(\L/\R)$ equals the determinant of $C$. To each root system we
associate two lists of integers: the degrees $d_1,\ldots,d_n$ and the
marks $c_1,\ldots,c_n$. Let $c$ denote the least common multiple of
the marks, and let $h$ denote the largest degree (called the Coxeter
number). Crystallographic root systems were classified by Cartan and
Killing in the late 1800s. Table~\ref{tab:weyl} lists the relevant
data for each Cartan--Killing type.

\begin{table}[h]
\begin{equation*}
\begin{array}{c|c|c|c|c|c}
\downstrut
\text{type} & \text{degrees} & \text{marks} & c & f & h \\
\hline
\upstrut
A_n & 2,3,\ldots,n+1 & 1,1,\ldots,1 &1& n+1 & n+1 \\
B_n \text{ and } C_n & 2,4,\ldots,2n &1,2,\ldots,2 &2& 2 & 2n \\
D_n & n, 2,4,\ldots,2(n-1) & 1,1,1,2,\ldots,2 &2& 4 & 2(n-1) \\
E_6 & 2,5,6,8,9,12 & 1,1,2,2,2,3  &6&3 & 12\\
E_7 & 2,6,8,10,12,14,18 & 1,2,2,2,3,3,4 &12& 2  &18\\
E_8 & 2,8,12,14,18, 20, 24, 30 & 2,2,3,3,4,4,5,6 &60& 1 & 30 \\
F_4 & 2,6,8,12 & 2,2,3,4 &12& 1 & 12\\
G_2 & 2,6 & 2,3 &6& 1 & 6\\
\end{array}
\end{equation*}
\caption{Root system data.}
\label{tab:weyl}
\end{table}

The fundamental alcove $\Delta$ is the simplex whose vertices
are $0$ and $\omega_i/c_i$ for all $1\leq i\leq n$. Suter
\cite{suter} observed (using an argument identical to the one in
Section~\ref{sec:qbinomial}) that the Ehrhart series of $\Delta$
has the form
\begin{equation*}
\sum_{b\geq 0} t^b \#(\L\cap b\Delta) = \frac{1}{(1-t)(1-t^{c_1})\cdots
(1-t^{c_n})}.
\end{equation*}
We associate to each root system a Weyl group $W$, whose size is
the product of the degrees: $\#W=d_1\cdots d_n$. If $\gcd(b,c)=1$
then Suter used known properties of the degrees and marks to observe that
\begin{equation*}
\#(\L\cap b\Delta) = \frac{f}{\#W} \prod_{i=1}^n (b+d_i-1)
\quad\text{when $\gcd(b,c)=1$},
\end{equation*}
and he asked for a conceptual explanation. Haiman \cite[Theorem~7.4.4]{haiman} had earlier proved a similar formula for the intersection of $b\Delta$ with the
root lattice $\R$. Recall that the vertices of $\Delta$ are $0$ and
$\omega_i/c_i$ for $1\leq i\leq n$. Since $c$ is the least common
multiple of the $c_i$, we observe that the vertices of $c\Delta$ are
in the weight lattice $\L$. Then since $f=\#(\L/\R)$, the vertices
of $cf\Delta$ are in the root lattice $\R$. Haiman used Burnside's
counting lemma, the Shephard--Todd formula for the dimensions of fixed
spaces, and Dirichlet's theorem on primes
in arithmetic progressions to prove that\footnote{Haiman observed
case by case that $\gcd(b,h)=1$ implies $\gcd(b,cf)=1$, where $h$
is the Coxeter number. Thiel \cite[Lemma 8.2]{thiel} gave a uniform
explanation of this phenomenon by applying a result of Bessis on the
degree of the generalized Lyashko--Looijenga morphism.}
\begin{equation*}
\#(\R\cap b\Delta) = \frac{1}{\#W} \prod_{i=1}^n (b+d_i-1)\quad
\text{when $\gcd(b,cf)=1$.}
\end{equation*}
Since $\gcd(b,cf)=1$ implies $\gcd(b,c)=1$, combining Suter's and
Haiman's results gives
\begin{equation*}
\#(\L\cap b\Delta) = f \cdot \#(\R\cap b\Delta) \quad\text{when
$\gcd(b,cf)=1$.}
\end{equation*}
This formula looks reasonable because $f=\#(\L/\R)$. One can imagine a more direct proof along the lines of the argument in Section~\ref{sec:rational_cat}. Humphreys 
\cite[Section~4.5]{humphreys} showed that there is a group of size $f$ acting by
symmetries on the fundamental alcove $\Delta$ and hence on the set
$\L\cap b\Delta$. If $\gcd(b,cf)=1$ then presumably each orbit of
this action has size $f$ and contains a unique point of $\R$. This would provide a conceptual explanation for Suter's formula.

Since $\#W=d_1\cdots d_n$, Haiman's formula can also be written as
\begin{equation*}
\#(\R\cap b\Delta) = \prod_{i=1}^n \frac{b-1+d_i}{d_i} \quad \text{
when $\gcd(b,cf)=1$}.
\end{equation*}
In type $A_{a-1}$ we have $c=1$, $f=a$ and
$(d_1,\ldots,d_{a-1})=(2,3,\ldots,a)$, so this formula becomes
\begin{equation*}
\#(\R\cap b\Delta) =
\frac{(b-1+2)(b-1+3)\cdots(b-1+a)}{(2)(3)\cdots(a)} = \frac{1}{a}
\binom{a-1+b}{a-1},
\end{equation*}
which is the rational Catalan number
$\Cat(a,b)$. Based on this analogy it makes sense to define the
following ``rational $q$-Catalan number'' for root systems:
\begin{equation*}
\Cat(W,b)_q:=  \prod_{i=1}^n \frac{[b-1+d_i]_q}{[d_i]_q} \quad \text{
when $\gcd(b,cf)=1$}.
\end{equation*}
We don't know where this formula first appeared, but the denominator
$[d_1]_q\cdots [d_n]_q$ is a well-known $q$-analogue of the Weyl
group cardinality coming from invariant theory. Gordon and Griffeth \cite{gordon} used the representation
theory of rational Cherednik algebras to prove that $\Cat(W,b)_q$ is
in $\NN[q]$. Galashin, Lam, Trinh and Williams~\cite{gltw} gave a new
proof that $\Cat(W,b)_q\in\ZZ[q]$ using Hecke algebras and generalized
Kazhdan--Lusztig polynomials, though their proof does not show that
the coefficients are positive. Neither proof gives a combinatorial
interpretation of the coefficients.\footnote{Both \cite{gordon}
and \cite{gltw} define $\Cat(W,b)_q$ in greater generality for
noncrystallographic reflection groups, where the definition is
more complicated.}

Based on the results in this paper, it is reasonable to conjecture
that there exists a Johnson statistic $\Jsf:\R\to\ZZ$ on the root
lattice $\R$ with the property that
\begin{equation*}
\Cat(W,b)_q= \sum_{\x\in (\R\cap b\Delta)} q^{\Jsf(\x)}  \quad \text{
when $\gcd(b,cf)=1$}.
\end{equation*}
In order to construct $\Jsf$ we should begin with a linear
``tilted height'' function $\Tsf:\L\to\ZZ$ satisfying
$\Tsf(\omega_i/c_i)=w_i\in\ZZ$ for some integers $w_1,\ldots,w_n$,
so that
\begin{equation*}
\Tsf(x_1\omega_1+\cdots+x_n\omega_n) = x_1w_1c_1+\cdots+x_nw_nc_n.
\end{equation*}
If $(\L\cap b\Delta)_q$ is the tilted height enumerator then the
weighted Ehrhart series is
\begin{equation*}
\sum_{b\geq 0} t^b (\L\cap b\Delta)_q =
\frac{1}{(1-t)(1-(q^{w_1}t)^{c_1})\cdots (1-(q^{w_n}t)^{c_n})}.
\end{equation*}
When $\gcd(b,cf)=1$ the formula $\#(\L\cap b\Delta) = f \cdot \#(\R\cap
b\Delta)$ should generalize to
\begin{equation*}
(\L\cap b\Delta)_q = f(q) \cdot \Cat(W,b)_q
\end{equation*}
for some appropriate $q$-analogue $f(q)\in\ZZ[q]$ of the index
of connection. Hopefully the statistic~$\Jsf$ will be a small
modification of $\Tsf$.

We demonstrate the feasibility of this approach in types $B_2$
and $G_2$.

\begin{figure}[h]
\begin{center}
\includegraphics[scale=0.9]{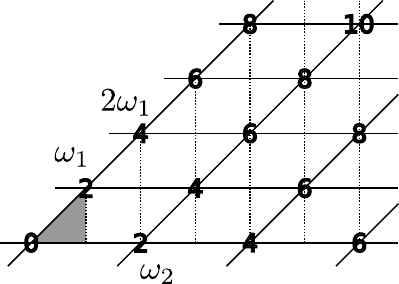}
\end{center}
\caption{Johnson statistic of type $B_2$.}
\label{fig:johnson_b2}
\end{figure}

In type $B_2$ we have $(d_1,d_2)=(2,4)$, $(c_1,c_2)=(1,2)$, $c=2$
and $f=2$. Let $\omega_1$, $\omega_2$ be the fundamental basis of the
weight lattice $\L$, so the vertices of the fundamental alcove $\Delta$
are $0$, $\omega_1$, $\omega_2/2$. The simple roots $\alpha_1$, $\alpha_2$
(expressed in weight coordinates) are the columns of the Cartan matrix
\begin{equation*}
C=\begin{pmatrix} \phantom{-}2 &-2 \\ -1 &\phantom{-}2\end{pmatrix},
\end{equation*}
that is, $\alpha_1=2\omega_1-\omega_2$ and
$\alpha_2=-2\omega_1+2\omega_2$. Let us take $(w_1,w_2)=(2,1)$ so the
tilted height function is $\Tsf(x_1\omega_1+x_2\omega_2)=2x_1+2x_2$
and the tilted height generating function is
\begin{equation*}
\sum_{b\geq 0} t^b (\L\cap b\Delta)_q =
\frac{1}{(1-t)(1-(q^2t)^1)(1-(q^1t)^2)}.
\end{equation*}
Since $cf=4$, the rational $q$-Catalan numbers are
\begin{equation*}
\Cat(B_2,b)_q = \frac{[b-1+2]_q[b-1+4]_q}{[2]_q[4]_q} \quad\text{when
$\gcd(b,4)=1$.}
\end{equation*}
One can check that
\begin{equation*}
\Cat(B_2,b)_q = \sum_{\x\in (\R\cap b\Delta)} q^{\Tsf(\x)} \quad
\text{when $\gcd(4,b)=1$},
\end{equation*}
so in this case the tilted height statistic is itself a Johnson
statistic. We also have
\begin{equation*}
\sum_{\x\in\L\cap b\Delta} q^{\Tsf(\x)} = (1+q^2) \cdot \sum_{\x\in
\R\cap b\Delta} q^{\Tsf(\x)}\quad \text{when $\gcd(4,b)=1$},
\end{equation*}
so the correct $q$-analogue of $f=2$ is $f(q)=1+q^2$. 
Figure~\ref{fig:johnson_b2} shows the intersection of $\L$ with the positive
cone, labeled by tilted height. Note that the root lattice has index
$2$ and is generated by the points $2\omega_1=2\alpha_1+\alpha_2$
and $\omega_2=\alpha_1+\alpha_2$.

\begin{figure}[h]
\begin{center}
\includegraphics[scale=0.9]{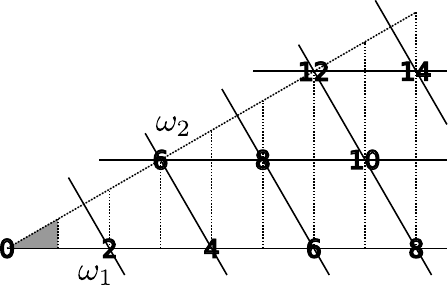}
\end{center}
\caption{Johnson statistic of type $G_2$.}
\label{fig:johnson_g2}
\end{figure}

In type $G_2$ we have $(d_1,d_2)=(2,6)$, $(c_1,c_2)=(2,3)$, $c=6$
and $f=1$. Let $\omega_1$, $\omega_2$ be the fundamental basis of the
weight lattice $\L$, so the vertices of the fundamental alcove $\Delta$
are $0$, $ \omega_1/2$, $\omega_2/3$. The simple roots $\alpha_1$, $\alpha_2$
(expressed in weight coordinates) are the columns of the Cartan matrix
\begin{equation*}
C=\begin{pmatrix} \phantom{-}2 &-1 \\ -3 &\phantom{-}2\end{pmatrix},
\end{equation*}
that is, $\alpha_1=2\omega_1-3\omega_2$ and
$\alpha_2=-3\omega_1+2\omega_2$. Let us define $(w_1,w_2)=(1,2)$ so
the tilted height function is $\Tsf(x_1\omega_1+x_2\omega_2)=2x_1+6x_2$
and the tilted height generating function is
\begin{equation*}
\sum_{b\geq 0} t^b (\L\cap b\Delta)_q =
\frac{1}{(1-t)(1-(q^1t)^2)(1-(q^2t)^3)}.
\end{equation*}
Since $cf=6$, the rational $q$-Catalan numbers are
\begin{equation*}
\Cat(G_2,b)_q = \frac{[b-1+2]_q[b-1+6]_q}{[2]_q[6]_q} \quad\text{when
$\gcd(6,b)=1$.}
\end{equation*}
One can check that
\begin{equation*}
\Cat(G_2,b)_q = \sum_{\x\in (\R\cap b\Delta)} q^{\Tsf(\x)} \quad
\text{when $\gcd(6,b)=1$},
\end{equation*}
so, again, the tilted height statistic is itself a Johnson
statistic. In this case the index of connection is $f=1$ so the root
and weight lattices coincide. Figure~\ref{fig:johnson_g2} shows the
intersection of~$\L$ with the positive cone, labeled by tilted height.

The problem of Johnson statistics is likely to be easier 
in types beyond $A$ because the index of connection is bounded above
by $4$.

\section*{Acknowledgements} We thank the following people for helpful input: Francois Bergeron, Benjamin Braun, Theresia Eisenk\"olbl, Eugene Gorsky, Christian Krattenthaler, Daniel Provder, Michael Schlosser, Nathan Williams, and the anonymous referees. We thank the anonymous referees and ChatGPT for help with the revised version of the paper. We thank Heather Armstrong for help with LaTeX. This project was partially supported by the SPP 2458 ``Combinatorial Synergies'', funded by the Deutsche Forschungsgemeinschaft (DFG, German Research Foundation).

\bibliographystyle{amsplain}
\bibliography{refs}

\end{document}